\documentclass{amsart}
\usepackage{amsmath}
\usepackage{amssymb}
\usepackage{amsthm} 
\usepackage{pxfonts}
\usepackage{times}
\usepackage{bbm}
\usepackage{graphicx}
\usepackage{algorithm}
\usepackage{algorithmic}
\usepackage{clrscode}
\usepackage{fullpage}

\newtheorem{cor}{Corollary}
\newtheorem{thm}{Theorem}
\newtheorem{lem}{Lemma}

\begin{document}
\pagestyle{plain}
\newenvironment{frcseries}{\fontfamily{frc} \selectfont}{}
\newcommand{\textfrc}[1]{{\frcseries #1}}
\newcommand{\mathfrc}[1]{\text{\textfrc{#1}}}

\title{Improved Approximation Guarantees for Sublinear-Time Fourier Algorithms}
%
\author{M. A. Iwen \\
Duke University, Box 90320 \\
Durham, NC 27708-0320\\
Email:  markiwen@math.duke.edu}

\maketitle
\thispagestyle{empty}

\begin{abstract}
In this paper modified variants of the sparse Fourier transform algorithms from \cite{JournalDSFT} are presented which improve on the approximation error bounds of the original algorithms.  In addition, simple methods for extending the improved sparse Fourier transforms to higher dimensional settings are developed.  As a consequence, approximate Fourier transforms are obtained which will identify a near-optimal $k$-term Fourier series for any given input function, $f: [0, 2 \pi]^D \rightarrow \mathbbm{C}$, in $O \left(k^2 \cdot D^4 \right)$ time (neglecting logarithmic factors).  Faster randomized Fourier algorithm variants with runtime complexities that scale linearly in the sparsity parameter $k$ are also presented.
\end{abstract}

\section{Introduction}
\label{sec:intro}

This paper develops fast methods for finding near-optimal nonlinear approximations to the Fourier transform of a given function $f: [0, 2\pi]^D \rightarrow \mathbbm{C}$.  Suppose that $f$ is a bandlimited function so that $\hat{f} \in \mathbbm{C}^{N^D}$, where $N^D$ is large.  An optimal $k$-term trigonometric approximation to $f$ is given by 
\begin{equation}
f^{\rm opt}_k (\vec{x}) = \sum^k_{j = 1} \hat{f} \left( \vec{\omega}_j \right) \mathbbm{e}^{\mathbbm{i} \vec{\omega}_j \cdot \vec{x} }
\label{BestkTerm}
\end{equation}
where $\vec{\omega}_1, \dots, \vec{\omega}_{N^D} \in [1-N/2, N/2]^D \cap \mathbbm{Z}^D$ are ordered by the magnitudes of their Fourier coefficients so that
$$\big| \hat{f}(\vec{\omega}_1) \big| \geq \big| \hat{f}( \vec{\omega}_2 ) \big| \geq \dots \geq \big| \hat{f}(\vec{\omega}_{N^D}) \big|.$$
The optimal $k$-term approximation error is then $\| f - f^{\rm opt}_k \|_2 = \| \hat{f} - \hat{f}^{\rm opt}_k \|_2$.  Suppose $k \in \mathbbm{N}$ is given.  The goal of this paper is to develop Fourier approximation schemes that are guaranteed to always return a near-optimal trigonometric polynomial, $y_k: [0, 2 \pi]^D \rightarrow \mathbbm{C}^{N^D}$, having
$\| f - y_k \|_2 \approx \| f - f^{\rm opt}_k \|_2$.  Furthermore, we require that the developed schemes are fast, with runtime complexities that scale polylogarithmically in $N^D$ and at most quadratically in $k$.  Such Fourier algorithms will then be able to accurately approximate the Fourier transform of an arbitrarily given function (i.e., with no a priori assumptions regarding ``smoothness'') much more quickly than a standard Fast Fourier Transform (FFT) methods \cite{FFT,BoydAl} whenever $N^D >> k$ is large.  More specifically, the developed schemes will lead to Fourier approximation algorithms with runtime complexities that scale \textit{polynomially} in $D$, as opposed to \textit{exponentially}.

The Fourier approximation techniques developed in this paper are improvements of the techniques introduced in \cite{JournalDSFT}.  As an example, suppose for simplicity that $f: [0, 2 \pi] \rightarrow \mathbbm{C}$ is a bandlimited function of only one variable so that $\hat{f} \in \mathbbm{C}^N$.  Furthermore, let $k < N$ be given.  The main theorem in \cite{JournalDSFT} implicitly proves that $O(k^2 \log^4 N)$ function evaluations and runtime are sufficient to produce a sparse approximation, $\hat{y}_k$, to $\hat{f}$ satisfying 
$$\left\| \hat{f} - \hat{y}_k \right\|_2 \leq \left\| \hat{f} - \hat{f}^{\rm opt}_k \right\|_2 + 3\sqrt{k} \left\| \hat{f} - \hat{f}^{\rm opt}_k \right\|_1,$$
where $f^{\rm opt}_k$ is defined as in Equation~\ref{BestkTerm}.  This error bound is unsatisfying for several reasons.  Principally, if many of the Fourier coefficients of $f$ are roughly the same magnitude the approximation error above can actually increase with $k$, the number of nonzero terms in the sparse approximation $\hat{y}_k$.  If nothing else, we would like to improve these error guarantees so that additional computational effort can always be counted on to yield better sparse Fourier approximations.

Let $p, q \in [1, \infty)$.  We will say that $\vec{y} \in \mathbbm{C}^N$ satisfies an $l^p$, $l^q / k^{1/q - 1/p}$ error bound with respect to $\hat{f} \in \mathbbm{C}^N$ if
\begin{equation}
\left\| \hat{f} - \vec{y} \right\|_p \leq \left\| \hat{f} - \hat{f}^{\rm opt}_k \right\|_p + \frac{\left\| \hat{f} - \hat{f}^{\rm opt}_k \right\|_q}{k^{1/q - 1/p}}.
\label{eqn:lplqError}
\end{equation}
More generally, we will refer to any error bound of the form given in Equation~\ref{eqn:lplqError} as an \textit{instance optimal error bound} for $\hat{f}$.  In this paper the result discussed in the previous paragraph is improved by showing that $O(k^2 \log^4 N)$ function samples and runtime are sufficient to produce a sparse approximation satisfying an $l^2$,~$l^1/\sqrt{k}$ error bound with respect to the Fourier transform of any $N$-bandlimited function $f: [0, 2 \pi] \rightarrow \mathbbm{C}$.  This decreases the ``$\sqrt{k} \left\| \hat{f} - \hat{f}^{\rm opt}_k \right\|_1$'' term in the previous error bound \cite{JournalDSFT} by a multiplicative factor of $k$.  Furthermore, faster randomized methods are also presented which are capable of achieving the same type of approximation errors typically achieved by slower algorithms based on the restricted isometry property \cite{NearOpt, HolgerCSBook} with high probability, despite utilizing a similar number of function samples.

\subsection{Results and Related Work}

Over the past few years, results concerning matrices with the Restricted Isometry Property (RIP) have allowed methods to be developed which can accurately approximate the Fourier transform of a function despite being given access to only a very small number of samples.  Informally, an $m \times N$ matrix $\mathcal{M}$ has the RIP of order $k \in \mathbbm{N}$ if it acts as a near isometry for all vectors, $\vec{x} \in \mathbbm{C}^N$, which contain at most $k$ nonzero entries.  Particularly important for our purposes is that RIP matrices of order $2k$ serve as good measurement matrices for sparsely approximating vectors in $\mathbbm{C}^N$. Suppose $\mathcal{M}$ is an $m \times N$ matrix with the RIP of order $2k$.  Then, for any $\vec{x} \in \mathbbm{C}^N$, a variety of computational methods including $l^1$-minimization \cite{CS1, CS4, NearOpt}, Orthogonal Matching Pursuit \cite{CS2,OMPvsBP}, Regularized Orthogonal Matching Pursuit \cite{ROMP,ROMPstable}, Iterative Hard Thresholding \cite{HardThreshforCS}, etc., will take $\mathcal{M} \vec{x}$ as input and subsequently output another vector, $\vec{y} \in \mathbbm{C}^N$, satisfying an instance optimal error bound with respect to $\vec{x}$ (e.g., an $l^2, l^1 / \sqrt{k}$ error bound).  Hence, any linear operator satisfying an appropriate RIP condition can serve as an efficient measurement operator capable of capturing sufficient information about any input vector in order to allow it to be accurately approximated.

The most pertinent RIP result to approximate Fourier recovery as considered here states that a rectangular matrix constructed by randomly selecting a small set of rows from an $N \times N$ inverse discrete Fourier transform matrix will have the RIP with high probability.  The following theorem was proven in \cite{CSb1} and subsequently generalized and improved in \cite{HolgerStableRecov}. 
\begin{thm} \textrm{\textbf{(See \cite{CSb1}).}}
Suppose we select $m$ rows uniformly at random from the rescaled $N \times N$ Inverse Discrete Fourier Transform (IDFT) matrix $\frac{1}{\sqrt{m}}\Psi^{-1}$, where
$$\left( \Psi^{-1} \right)_{i,j} = \frac{e^{\frac{2 \pi \mathbbm{i} \cdot i \cdot j}{N}}}{\sqrt{N}},$$
and form the $m \times N$ submatrix $\mathcal{M}$.  If $m$ is $\Omega \left(k \cdot \log N \cdot \log^2 k \cdot \log(k \log N) \right)$ then $\mathcal{M}$ will have the RIP of order $k$ with high probability.
\label{thm:FourierRIPrand}
\end{thm} 

Let $\Psi$ be the $N \times N$ Discrete Fourier Transform (DFT) matrix defined by $\Psi_{i,j} = \frac{1}{N} \cdot e^{\frac{2 \pi \mathbbm{i} \cdot i \cdot j}{N}}$, $f: [0, 2 \pi] \rightarrow \mathbbm{C}$ be a given function, and $\vec{f} \in \mathbbm{C}^N$ be the vector of $N$ equally spaced samples from $f$ on $[0, 2 \pi]$.  In this case Theorem~\ref{thm:FourierRIPrand} tells us that collecting the $m$ function samples determined by $\mathcal{M} \Psi \vec{f}$ will be sufficient to accurately approximate the discrete Fourier transform of $\vec{f}$ with high probability.  More precisely, if $\mathcal{M} \hat{\vec{f}} = \mathcal{M} \Psi \vec{f}$ is input to a recovery algorithm known as CoSaMP \cite{COSAMP} the following theorem holds.
\begin{thm} \textrm{\textbf{(See \cite{COSAMP}).}}
Suppose that $\mathcal{M}$ is a $m \times N$ measurement matrix formed by selecting $m = \Theta(k \cdot \log^4 N)$ rows from the $N \times N$ IDFT matrix, $\Psi^{-1}$, uniformly at random.  Furthermore, assume that $\mathcal{M}$ satisfies the RIP of order $2k$\footnote{Note that this is true with high probability by Theorem~\ref{thm:FourierRIPrand}.}.  Fix precision parameter $\eta \in \mathbbm{R}$ and let $\vec{U} = \mathcal{M} \Psi \vec{f}$ be measurements collected for any given $\vec{f} \in \mathbbm{C}^N$.  Then, when executed with $\vec{U}$ as input, CoSaMP will output a $2k$-sparse vector, $\vec{y} \in \mathbbm{C}^N$, satisfying 
$$\left\| \hat{\vec{f}} - \vec{y} \right\|_2 \leq {\rm Const} \cdot \max \left\{ \eta, \frac{1}{\sqrt{k}} \cdot \left\| \hat{\vec{f}} - \hat{\vec{f}}^{\rm ~opt}_k \right\|_1 \right\},$$
where $\hat{\vec{f}}^{\rm ~opt}_k$ is a best possible $k$-term approximation for $\hat{\vec{f}} = \Psi \vec{f}$.  The required runtime is $O \left(N \log N \cdot \log \left( \left\| \vec{f}~ \right\|_2 / \eta \right) \right)$.
\label{thm:Cosamp}
\end{thm}

In effect, Theorem~\ref{thm:Cosamp} promises that CoSaMP will locate $2k$ of the dominant entries in $\hat{\vec{f}}$ if given access to $\Theta(k \cdot \log^4 N)$ samples from $f$.  If $\hat{\vec{f}}$ contains $2k$ significant frequencies whose Fourier coefficients collectively dominate all others combined, then these most significant frequencies will be found and their Fourier coefficients will be well approximated.  If $\hat{\vec{f}}$ has no dominant set of $2k$ entries then CoSaMP will return a sparse representation which is guaranteed only to be trivially bounded.
However, in such cases sparse Fourier approximation is a generally hopeless task anyways and a bounded, albeit poor, sparse representation is the best one can expect.  In any case, as long as the random function samples correspond to a matrix with the RIP, CoSaMP will output a vector satisfying an instance optimal error bound with respect to $\hat{\vec{f}}$.  However, the required runtime will always be $\Omega(N)$.  More generally, all existing Fourier recovery methods based on RIP conditions have superlinear runtime complexity in $N$.

Other existing Fourier algorithms for approximating $\hat{\vec{f}} \in \mathbbm{C}^N$ given sampling access to $\vec{f} \in \mathbbm{C}^N$ work by utilizing random sampling techniques \cite{AAFFT1,AAFFT2}. These approaches simultaneously obtain both instance optimal error guarantees, and runtime complexities that scale sublinearly in $N$. However, they generally also require more function samples than recovery algorithms which utilize matrices satisfying the RIP.  A variant of the following Fourier sampling theorem, concerning the sparse approximation of $\hat{\vec{f}}$ provided sampling access to $f: [0, 2 \pi] \rightarrow \mathbbm{C}$, is proven in \cite{AAFFT2}.
\begin{thm} \textrm{\textbf{(See \cite{AAFFT2}).}}
Fix precision parameters $\eta,\tau \in \mathbbm{R}^{+}$ and probability parameter $\lambda \in (0,1)$.  There exists a randomized sampling algorithm which, when given sampling access to an input signal $\vec{f} \in \mathbbm{C}^N$, outputs a $k$-sparse representation $\vec{y}$ for $\hat{\vec{f}}$ satisfying 
$$\left\| \hat{\vec{f}} - \vec{y} \right\|_2 \leq \sqrt{1 + \tau} \cdot \max \left\{ \eta, \left\| \hat{\vec{f}} - \hat{\vec{f}}^{\rm ~opt}_k \right\|_2 \right\}$$
with probability at least $1 - \lambda$.  Here $\hat{\vec{f}}^{\rm ~opt}_k$ is a best possible $k$-sparse representation for $\hat{\vec{f}}$.  Both the runtime and sampling complexities are bounded above by 
$$k \cdot \left( \log \left( \frac{1}{\lambda} \right), \log \left( \frac{1}{\eta} \right), \log  \| \textbf{A} \|_2, \log N ,\frac{1}{\tau} \right)^{O(1)}.$$  
\label{thm:AAFFT}
\end{thm}
It is important to note that the probabilistic guarantee of recovering an accurate sparse representation provided by Theorem~\ref{thm:AAFFT} is a nonuniform \textit{per signal} guarantee.  In contrast, Fourier approximation procedures which rely on RIP matrices provide uniform probability guarantees \textit{for all} possible input vectors.  If a set of sample positions corresponds to an $N \times N$ IDFT submatrix with the RIP property, those sample positions will allow the accurate Fourier approximation of all possible input vectors $\vec{f} \in \mathbbm{C}^N$. 

\begin{table} [tb]
\begin{center}
\begin{tabular}{|c|c|c|c|c|}
	\hline
\textbf{Fourier Result} & \textbf{w.h.p./D} & \textbf{Runtime} & \textbf{Function Samples} & \textbf{Error Guarantee} \\
\hline
Theorem~\ref{thm:Recon1} & D & $O(N \cdot k \cdot \log^2 N )$ & $O(k^2 \cdot \log^2 N)$ & $l^2$,~$l^1/\sqrt{k}$ \\
\hline
CoSaMP \cite{COSAMP} & $\approx$D & $O(N \cdot \log N)$ & $O(k \cdot \log^4 N)$ \cite{CSb1, HolgerCSBook} & $l^2$,~$l^1/\sqrt{k} + \eta$ \\
\hline
Corollary~\ref{cor:RandRecon1} & w.h.p. & $O \left(N \cdot \log N \right)$ & $O \left(k \cdot \log^2 N \right)$ & $l^2$,~$l^1/\sqrt{k}$ \\
\hline
\hline
Theorem~\ref{thm:Recon2} & D & $O(k^2 \cdot \log^4 N)$ & $O(k^2 \cdot \log^4 N)$ & $l^2$,~$l^1/\sqrt{k}$ \\
\hline
Sparse Fourier \cite{AAFFT2} & w.h.p. & $O \left(k \cdot \log^{O(1)}(N) \right)$ & $O \left(k \cdot \log^{O(1)}(N) \right)$ & $l^2$,~$l^2 + \eta$ \\
\hline
Corollary~\ref{cor:Recon2} & w.h.p. & $O\left( k \cdot \log^5 N \right)$ & $O \left(k \cdot \log^4 N \right)$ & $l^2$,~$l^1/\sqrt{k}$ \\
\hline
\hline
Optimal Algorithm & D & $\Omega \left( k \right)$ & $\Omega \left( k \right)$ \cite{BestkTerm} & $l^2$,~$l^1/\sqrt{k}$\\
\hline
\end{tabular}
\end{center}
\caption{Sparse Fourier Approximation Algorithms with Robust Recovery Guarantees}
\label{tab:Compare1}
\end{table}

In this paper several Fourier algorithms are developed which obtain instance optimal approximation guarantees while also improving on various aspects of the previously mentioned approaches.  See Table~\ref{tab:Compare1} for a comparison of the results obtained herein with Theorems~\ref{thm:Cosamp} and~\ref{thm:AAFFT} when applied to the problem of approximating the Fourier transform, $\hat{f} \in \mathbbm{C}^N$, of an $N$-bandwidth function $f:  [0, 2 \pi] \rightarrow \mathbbm{C}$.  The first column of Table~\ref{tab:Compare1} lists the Fourier results considered, while the second column lists whether the recovery algorithm in question guarantees an instance optimal output Deterministically (D), or With High Probability (w.h.p.) per signal.  Note that CoSaMP\footnote{We used CoSaMP as a representative for all RIP based recovery algorithms because, for the purposes of Table~\ref{tab:Compare1} at least, it matches the currently best achievable runtime, sampling, and error bound performance characteristics of all the other previously mentioned RIP-based methods in the Fourier setting.} has an ``$\approx D$'' listed in its second column.  This denotes that the RIP results utilized in Theorem~\ref{thm:Cosamp} provide a uniform probability guarantee, although no explicit constructions of RIP matrices satisfying these bounds are currently known.  The third and fourth columns of Table~\ref{tab:Compare1} contain the sampling and runtime complexities of the algorithms, respectively.  For simplicity some of the bounds were simplified by ignoring precision parameters, etc.\footnote{The $O(N \cdot \log N)$ runtime listed for Corollary~\ref{cor:RandRecon1} will hold if $k$ is $O(N / \log^2 N)$.  More generally, the runtime will always be $O(N \cdot \log^3 N)$.}.  Finally, the fifth column of Table~\ref{tab:Compare1} lists the instance optimal approximation guarantees achievable by each algorithm when budgeted the number of samples and time listed in the third and fourth columns.  The ``$+ \eta$'' in the CoSaMP and Sparse Fourier rows remind us that their error bounds are good up to an additive precision parameter.

The last row of Table~\ref{tab:Compare1} lists lower bounds for the runtime and sampling complexity of any algorithm guaranteed to achieve an instance optimal $l^2$,~$l^1/\sqrt{k}$ Fourier approximation error (see \cite{BestkTerm}).  Note that all six approaches have sampling complexities containing additional multiplicative logarithmic factors of $N$ beyond the stated lower sampling bound \footnote{Both the sampling and runtime complexities of the Sparse Fourier algorithm presented in \cite{AAFFT2} scale like $\Omega( k \cdot \log^5 N )$.}.  The lowest overall sampling complexity is achieved by Corollary~\ref{cor:RandRecon1}, although, it is achieved at the expense of a weak nonuniform ``w.h.p.'' approximation probability guarantee.  Similarly, Corollary~\ref{cor:Recon2} improves on the previous sampling complexity of the sparse Fourier algorithm in \cite{AAFFT2} while at least matching its runtime complexity\footnote{It must be remembered, however, that the algorithm presented in \cite{AAFFT2} enjoys a stronger approximation error guarantee up to its additive/multiplicative precision parameters.}.  Finally, to the best of the author's knowledge, Theorem~\ref{thm:Recon2} obtains the best available runtime of any existing deterministic Fourier approximation algorithm which is guaranteed to achieve an instance optimal error guarantee.

The remainder of this paper is organized as follows:  In Section~\ref{sec:prelim} the notation utilized throughout the remainder of the paper is established.  Next, in Section~\ref{sec:makeM}, a number theoretic matrix construction is presented and analyzed.  Section~\ref{sec:FixedSig} explains how random submatrices of the presented number theoretic matrices can yield nonuniform probabilistic approximation guarantees, while Section~\ref{sec:Fourier} outlines a useful relationship between these matrices and the Fourier transform of a periodic function.  In Section~\ref{sec:FourierRecov} the matrices defined in Section~\ref{sec:makeM} are used to construct Fourier approximation algorithms with runtime complexities that scale superlinearly in $N$ (i.e., Theorem~\ref{thm:Recon1} and Corollary~\ref{cor:RandRecon1} are proven).  Next, in Section~\ref{sec:FourierRecov2}, the algorithms of Section~\ref{sec:FourierRecov} are modified into algorithms with runtime complexities that scale sublinearly in $N$ (i.e., Theorem~\ref{thm:Recon2} and Corollary~\ref{cor:Recon2} are proven).  In Section~\ref{sec:HighDim} a simple strategy is given for extending the results of the previous two sections to higher dimensional Fourier transforms.  Finally, a short conclusion is presented in Section~\ref{sec:Conc}.

\section{Notation and Setup}
\label{sec:prelim}

Below we will consider any function whose domain, $I$, is both ordered and countable to be a vector.  Let $\vec{x}:  I \rightarrow \mathbbm{C}$.  In this case we will say that $\vec{x} \in \mathbbm{C}^{|I|}$, and that $x_i = \vec{x}(i) \in \mathbbm{C}$ for all $i \in I$.  We will denote the $l^p$ norm of any such vector, $\vec{x}$, by
$$\| \vec{x} \|_p = \left( \sum_{i \in I} |x_i|^p \right)^{\frac{1}{p}},~\textrm{ for }p \in [1, \infty).$$
If $\vec{x}$ is an infinite vector (i.e., if $I$ is countably infinite), we will say that $\vec{x} \in l^p$ if $\| \vec{x} \|_p$ is finite.  Without loss of generality, we will assume that a given $\vec{x} \in \mathbbm{C}^N$ is indexed by $I = [0, N) \cap \mathbbm{Z}$ unless indicated otherwise.  The vector $\vec{\mathbbm{1}}_N \in \mathbbm{C}^N$ will always denote the vector of $N$ ones, and $\vec{0}_N \in \mathbbm{C}^N$ with always denote the vector of $N$ zeros.

For any given $\vec{x} \in \mathbbm{C}^{|I|}$ and subset $S \subseteq I$, we will let $\vec{x}_{S} \in \mathbbm{C}^{|I|}$ be equal to $\vec{x}$ on $S$ and be zero everywhere else.  Thus, 
$$\left( \vec{x}_{S} \right)_i = \left\{ \begin{array}{ll} x_i & \textrm{if } i \in S, \\ ~0 & {\rm otherwise} \end{array} \right..$$
Furthermore, for a given integer $k < |I|$, we will let $S^{\rm opt}_k \subset I$ be the first $k$ element subset of $I$ in lexicographical order with the property that $|x_s| \geq |x_t|$ for all $s \in S^{\rm opt}_k$ and $t \in I - S^{\rm opt}_k$.  Thus, $S^{\rm opt}_k$ contains the indexes of $k$ of the largest magnitude entries in $\vec{x}$.  Finally, we will define $\vec{x}^{\rm opt}_k$ to be $\vec{x}_{S^{\rm opt}_k}$, a best $k$-term approximation to $\vec{x}$.

In this paper we will be considering methods for approximating the Fourier series of an arbitrarily given periodic function, $f: [0, 2\pi]^D \rightarrow \mathbbm{C}$.  Following  convention, we will denote the Fourier transform of $f$ by $\hat{f}: \mathbbm{Z}^D \rightarrow \mathbbm{C}$, where
$$\hat{f} \left( \vec{\omega} \right) = \frac{1}{\left(2 \pi \right)^D} \int_{\vec{x} \in [0, 2\pi]^D} \mathbbm{e}^{- \mathbbm{i} \vec{\omega} \cdot \vec{x}} ~f(\vec{x})~ d \vec{x} ~~\textrm{ for all } \vec{\omega} \in \mathbbm{Z}^D.$$  
Note that $\hat{f}$ can be considered an infinite vector indexed by $\mathbbm{Z}^D$.  We also have the inverse relationship
$$f(\vec{x}) = \sum_{\vec{\omega} \in \mathbbm{Z}^D} \hat{f}\left( \vec{\omega} \right) ~\mathbbm{e}^{\mathbbm{i} \vec{\omega} \cdot \vec{x}} ~~\textrm{ for all } \vec{x} \in [0, 2\pi]^D.$$  
Thus, we learn $f$ in the process of approximating its Fourier transform.

Call each $\vec{\omega} \in \mathbbm{Z}^D$ a \textit{Fourier mode} or \textit{frequency}, and $\hat{f} \left( \vec{\omega} \right)$ its corresponding \textit{Fourier coefficient}.  Ultimately, we will restrict our attention to the Fourier modes of $f$ inside some finite bandwidth.  We will do this by identifying, and then estimating the Fourier coefficients of, the most energetic Fourier modes in $\left( -\left\lceil \frac{N}{2} \right\rceil, \left\lfloor \frac{N}{2} \right\rfloor \right]^D \cap \mathbbm{Z}^D$ for a given bandwidth value $N \in \mathbbm{N}$.  Toward this end, define the vector $\vec{\hat{f}} \in \mathbbm{C}^{N^D}$ by
$$\hat{f}_{\vec{\omega}} ~=~ \hat{f} \left( \vec{\omega} \right) ~~\textrm{ for all } \vec{\omega} \in \left( -\left\lceil \frac{N}{2} \right\rceil, \left\lfloor \frac{N}{2} \right\rfloor \right]^D \cap \mathbbm{Z}^D.$$
Similarly, define $\bar{\hat{f}}: \mathbbm{Z}^D \rightarrow \mathbbm{C}$ to be the Fourier transform of the related optimal bandlimited approximation to $f$.  More precisely, let
$$\bar{\hat{f}} \left( \vec{\omega} \right) = \hat{f}_{\left( -\left\lceil \frac{N}{2} \right\rceil, \left\lfloor \frac{N}{2} \right\rfloor \right]^D \cap \mathbbm{Z}^D} = \left\{ \begin{array}{ll} \hat{f} \left( \vec{\omega} \right) & \textrm{if } \vec{\omega} \in \left( -\left\lceil \frac{N}{2} \right\rceil, \left\lfloor \frac{N}{2} \right\rfloor \right]^D \cap \mathbbm{Z}^D, \\ ~~0 & {\rm otherwise} \end{array} \right.$$
for all $\omega \in \mathbbm{Z}^D$.  We will approximate $\hat{f}$ by approximating $\vec{\hat{f}}$.  However, in order to do so we must first construct a special class of matrices.

\section{A Specialized Measurement Matrix Construction}
\label{sec:makeM}

We consider $m \times N$ measurement matrices, $\mathcal{M}_{s_1,K}$, constructed as follows.  Select $K$ pairwise relatively prime integers beginning with a given $s_1 \in \mathbbm{N}$ and denote them by
\begin{equation}
s_1 < \cdots < s_K.
\label{def:sj}
\end{equation}
Produce a row $r_{j,h}$, where $j \in [1,K] \cap \mathbbm{N}$ and $h \in [0,s_j) \cap \mathbbm{N}$, 
in $\mathcal{M}_{s_1,K}$ for each possible residue of each $s_j$ integer.  The $n^{\rm th}$ entry of each $r_{j,h}$ row, $n \in [0,N) \cap \mathbbm{N}$, is given by
\begin{align}
(r_{j,h})_{n} & = \delta \left( (n-h) \textrm{ mod } s_j \right) \nonumber \\ & = \left\{ \begin{array}{ll} 1 & \textrm{if } n \equiv h \textrm{ mod } s_j \\ 0 & {\rm otherwise} \end{array} \right. .
\label{eqn:Def_r}
\end{align}
We then set
\begin{equation}
\mathcal{M}_{s_1,K} = \left( \begin{array}{l} r_{1,0} \\ r_{1,1}\\ \vdots \\ r_{1,s_{1}-1} \\ \vdots \\ r_{K,s_{K}-1} \\
\end{array} \right).
\label{eqn:Def_M}
\end{equation}
The result is an $\left( m = \sum^K_{j = 1} s_j \right) \times N$ matrix with binary entries.  See Figure~\ref{fig:example} for an example measurement matrix.

\begin{figure}
$$\textbf{---------------------------------------------------------------------------}$$
$$\begin{array}{llllllllll}
\hspace{-7pt} \mathbf{n} \in \mathbf{[0,N)} \cap \mathbbm{N}&& \mathbf{0} & \mathbf{1} & \mathbf{2} & \mathbf{3} & \mathbf{4} & \mathbf{5} & \mathbf{6} & \dots \\ 
\end{array}$$
$$\begin{array}{l} \mathbf{n} \equiv \mathbf{0}~\mathbf{mod}~\mathbf{2} \\ \mathbf{n} \equiv \mathbf{1}~\mathbf{mod}~\mathbf{2} \\ \mathbf{n} \equiv \mathbf{0}~\mathbf{mod}~\mathbf{3} \\ \mathbf{n} \equiv \mathbf{1}~\mathbf{mod}~3 \\ \mathbf{n} \equiv \mathbf{2}~\mathbf{mod}~\mathbf{3} \\ ~~~~~~~~~~\vdots \\ \mathbf{n} \equiv \mathbf{1}~\mathbf{mod}~\mathbf{5} \\ ~~~~~~~~~~\vdots \end{array} 
\begin{array}{l}  \\  \\  \\  \\  \\  \\  \\  \end{array} 
\left( \begin{array}{llllllll} 
1 & 0 & 1 & 0 & 1 & 0 & 1 & \dots \\ 
0 & 1 & 0 & 1 & 0 & 1 & 0 & \dots \\ 
1 & 0 & 0 & 1 & 0 & 0 & 1 & \dots \\ 
0 & 1 & 0 & 0 & 1 & 0 & 0 & \dots \\ 
0 & 0 & 1 & 0 & 0 & 1 & 0 & \dots \\
&&& \vdots &&&& \\
0 & 1 & 0 & 0 & 0 & 0 & 1 & \dots \\ 
&&& \vdots &&&& \\
\end{array} \right)$$ 
\caption{An Example Matrix created using $s_{1} = 2$, $s_{2} = 3$, $s_{3} = 5$, \dots}
\label{fig:example}
$$\textbf{---------------------------------------------------------------------------}$$
\end{figure}

The matrices constructed above using relatively prime integers have many useful properties.  As we shall see later in Section~\ref{sec:FourierRecov}, these properties cumulatively allow the accurate recovery of Fourier sparse signals.
We require two additional definitions before we may continue.  Let $n \in [0,N) \cap \mathbbm{N}$.  We define $\mathcal{M}_{s_1,K,n}$ to be the $K \times N$ matrix created by selecting the $K$ rows of $\mathcal{M}_{s_1,K}$ with nonzero entries in the $n^{\rm th}$ column.  Furthermore, we define $\mathcal{M'}_{s_1,K,n}$ to be the $K \times (N-1)$ matrix created by deleting the $n^{\rm th}$ column of $\mathcal{M}_{s_1,K,n}$.  Thus, we have 
\begin{equation}
\mathcal{M}_{s_1,K,n} = \left( \begin{array}{l} r_{1,~n \textrm{~mod~} s_1} \\ 
r_{2,~n \textrm{~mod~} s_2}\\ \vdots \\ r_{K,~n \textrm{~mod~} s_K}\\ 
\end{array} \right)
\label{def:Msub}
\end{equation}
and
\begin{equation}
\mathcal{M'}_{s_1,K,n} = \left( \begin{array}{lllllll} (r_{1,n \textrm{~mod~} s_1})_0 & (r_{1,n \textrm{~mod~} s_1})_1 & \dots & (r_{1,n \textrm{~mod~} s_1})_{n-1} & (r_{1,n \textrm{~mod~} s_1})_{n+1} & \dots & (r_{1,n \textrm{~mod~} s_1})_{N-1} \\ 
(r_{2,n \textrm{~mod~} s_2})_0 & (r_{2,n \textrm{~mod~} s_2})_1 & \dots & (r_{2,n \textrm{~mod~} s_2})_{n-1} & (r_{2,n \textrm{~mod~} s_2})_{n+1} & \dots & (r_{2,n \textrm{~mod~} s_2})_{N-1}\\ 
&&&&\hspace{-5pt} \vdots&& \\ 
(r_{K,n \textrm{~mod~} s_K})_0 & (r_{K,n \textrm{~mod~} s_K})_1 & \dots & (r_{K,n \textrm{~mod~} s_K})_{n-1} & (r_{K,n \textrm{~mod~} s_K})_{n+1} & \dots & (r_{K,n \textrm{~mod~} s_K})_{N-1}\\ 
\end{array} \right).
\label{def:M'sub}
\end{equation}
We have the following two lemmas.

\begin{lem}
Let $n, \bar{k} \in [0,N) \cap \mathbbm{N}$ and $\vec{x} \in \mathbbm{C}^{N-1}$.  Then, at most $\bar{k} \left \lfloor \log_{s_1} N \right \rfloor$ of the $K$ entries of $\mathcal{M'}_{s_1,K,n} \cdot \vec{x}$ will have magnitude greater than or equal to $\| \vec{x} \|_1 / \bar{k}$.
\label{lem1}
\end{lem}

\noindent \textit{Proof:}\\  

We have that
$$\left| \left\{ j ~\bigg|~ \left| \left(\mathcal{M'}_{s_1,K,n} \cdot \vec{x} \right)_j \right| \geq \frac{\| \vec{x} \|_1 }{\bar{k}} \right\} \right| \leq \frac{\bar{k}}{\| \vec{x} \|_1} \left\| \mathcal{M'}_{s_1,K,n} \cdot \vec{x} \right\|_1 \leq \bar{k} \cdot \| \mathcal{M'}_{s_1,K,n} \|_1$$
by the Markov Inequality.  Focusing now on $\mathcal{M'}_{s_1,K,n}$ we can see that
\begin{equation}
\left \| \mathcal{M'}_{s_1,K,n} \right\|_1 = \max_{l \in [0,N-1) \cap \mathbbm{N}} \sum^{K}_{j = 1} \left| \left(\mathcal{M'}_{s_1,K,n} \right)_{j,l} \right| = \max_{l \in [0,N-1) \cap \mathbbm{N}} \sum^{K}_{j = 1} \delta \left( (n - l) \textrm{~mod~} s_j\right) \leq \left \lfloor \log_{s_1} N \right \rfloor
\label{eqn:l1Bound}
\end{equation}
by the Chinese Remainder Theorem (see~\cite{NumTheory}).  The result follows.~~$\Box$ \\

\begin{lem}
Let $n,\tilde{k} \in [0,N) \cap \mathbbm{N}$, $S \subset [0,N) \cap \mathbbm{N}$ with $|S| \leq \tilde{k}$, and $\vec{x} \in \mathbbm{C}^{N-1}$.  Then, $\mathcal{M'}_{s_1,K,n} \cdot \vec{x}$ and $\mathcal{M'}_{s_1,K,n} \cdot \left( \vec{x} - \vec{x}_{S} \right)$ will differ in at most $\tilde{k} \left \lfloor \log_{s_1} N \right \rfloor$ of their $K$ entries.  
\label{lem2}
\end{lem}

\noindent \textit{Proof:} \\

We have that
$$\left| \left\{ j ~\bigg|~ \left(\mathcal{M'}_{s_1,K,n} \cdot \vec{x} \right)_j  \neq \left( \mathcal{M'}_{s_1,K,n} \cdot \left( \vec{x} - \vec{x}_{S} \right) \right)_j \right\} \right| = \left| \left\{ j ~\bigg|~ \left(\mathcal{M'}_{s_1,K,n} \cdot \vec{x}_{S} \right)_j  \neq 0 \right\} \right| \leq \left| \left\{ j ~\bigg|~ \left(\mathcal{M'}_{s_1,K,n} \cdot \left(\vec{\mathbbm{1}}_{N-1}\right)_{S} \right)_j  \geq 1 \right\} \right|$$
since all the entries of $\mathcal{M'}_{s_1,K,n}$ are nonnegative integers.  Applying Lemma~\ref{lem1} with $\vec{x} = \left(\vec{\mathbbm{1}}_{N-1}\right)_{S}$ and $\bar{k} = \left \| \left(\vec{\mathbbm{1}}_{N-1}\right)_{S} \right \|_1 = |S|$ finishes the proof.~~$\Box$ \\

Combining these two Lemmas we obtain a general theorem concerning the accuracy with which we can approximate any entry of an arbitrary complex vector $\vec{x} \in \mathbbm{C}^N$ using only entries of $\mathcal{M}_{s_1,K} \cdot \vec{x}$.

\begin{thm}
Let $n,k,s_1 \in [0,N) \cap \mathbbm{N}$, $\epsilon^{-1} \in \mathbbm{N}^+ $, $c \in [2,\infty) \cap \mathbbm{N}$, and $\vec{x} \in \mathbbm{C}^{N}$.  Set $K = c \cdot (k / \epsilon) \left \lfloor \log_{s_1} N \right \rfloor + 1$.  Then, more than $\frac{c-2}{c} \cdot K$ of the $K$ entries of $\mathcal{M}_{s_1,K,n} \cdot \vec{x}$ will estimate $x_n$ to within $\frac{\epsilon \cdot \left\| \vec{x} - \vec{x}^{\rm opt}_{(k/\epsilon)} \right\|_1}{k}$ precision.
\label{thm:Mest}
\end{thm}

\noindent \textit{Proof:}\\  

Define $\vec{y} \in \mathbbm{C}^{N-1}$ to be $\vec{y} = \left( x_0, x_1, \dots, x_{n-1}, x_{n+1}, \dots, x_{N-1} \right)$.  We have that
$$\mathcal{M}_{s_1,K,n} \cdot \vec{x} ~=~ x_n \cdot \vec{\mathbbm{1}}_K + \mathcal{M'}_{s_1,K,n} \cdot \vec{y}.$$  Applying Lemma~\ref{lem2} with $\tilde{k} = (k/\epsilon)$ reveals that at most $(k/\epsilon) \left \lfloor \log_{s_1} N \right \rfloor$ entries of $\mathcal{M'}_{s_1,K,n} \cdot \vec{y}$ differ from $\mathcal{M'}_{s_1,K,n} \cdot \left( \vec{y} - \vec{y}^{\rm opt}_{(k / \epsilon)} \right)$.  Of the remaining $K - (k/\epsilon) \left \lfloor \log_{s_1} N \right \rfloor$ entries of $\mathcal{M'}_{s_1,K,n} \cdot \vec{y}$, at most $(k  / \epsilon) \left \lfloor \log_{s_1} N \right \rfloor$ will have magnitudes greater than or equal to $\epsilon \left\| \vec{y} - \vec{y}^{\rm opt}_{(k / \epsilon)} \right\|_1 / k$ by Lemma~\ref{lem1}.  Hence, at least 
$$K - 2(k/\epsilon) \left \lfloor \log_{s_1} N \right \rfloor \geq (c - 2) (k / \epsilon) \left \lfloor \log_{s_1} N \right \rfloor + 1 > \frac{c-2}{c} \cdot K$$
entries of $\mathcal{M'}_{s_1,K,n} \cdot \vec{y}$ will have a magnitude no greater than
$$\frac{\epsilon \cdot \left\| \vec{y} - \vec{y}^{\rm opt}_{(k / \epsilon)} \right\|_1}{k} \leq \frac{\epsilon \cdot \left\| \vec{x} - \vec{x}^{\rm opt}_{(k / \epsilon)} \right\|_1}{k}.$$
The result follows.~~$\Box$ \\

We will now study the number of rows, $m = \sum^{K}_{j=1} s_{j}$, in our measurement matrix under the Theorem~\ref{thm:Mest} assumption that $K = c \cdot (k / \epsilon) \left \lfloor \log_{s_1} N \right \rfloor + 1$ for some constant integer $c \in [2,\infty)$ and given values of $s_1 = (k/\epsilon),N \in \mathbbm{N}^+$.  Given this assumption concerning $K$, we wish to bound the smallest possible sum, $m$, resulting from all possible choices of pairwise relatively prime $s_j$ values.  We will do this by bounding $m$ for one particular set of $s_j$ values.

Let $p_{l}$ be the $l^{\rm th}$ prime natural number.  Thus, we have
\begin{equation}
p_{1} = 2, p_{2} = 3, p_{3} = 5, p_{4} = 7, \dots
\label{eqn:Primes}
\end{equation}
Next, define $q\in \mathbbm{N}$ so that
\begin{equation}
p_{q-1} < (k/\epsilon) \leq p_{q}.
\label{eqn:Def_q}
\end{equation}
We will use the first $K$ primes no smaller than $(k/\epsilon)$ to define our relatively prime $s_j$ values for the purposes of bounding $m$.  Hence, for the remainder of Section~\ref{sec:makeM} we will have
\begin{equation}
s_1 = \frac{k}{\epsilon} \leq p_{q} < s_2 = p_{q+1} < \dots < s_K = p_{q+K-1}.
\label{eqn:Def_k}
\end{equation}

It follows from results in \cite{IDSFA} that 
\begin{equation}
m = \sum^K_{j=1} s_j \leq \sum^{K-1}_{j=0} p_{q+j} = \frac{p^2_{q+K}}{2 \ln p_{q+K}} \cdot \left( 1 + O \left( \frac{1}{\ln p_{q+K}} \right) \right) -  \frac{p^2_{q}}{2 \ln p_{q}} \cdot \left( 1 + O \left( \frac{1}{\ln p_{q}} \right) \right).
\label{eqn:Sum}
\end{equation}
Furthermore, the Prime Number Theorem (see \cite{NumTheory}) tells us that
$$q = \frac{k}{\epsilon \cdot \ln (k/\epsilon)} \left( 1+ O\left( \frac{1}{\ln (k/\epsilon)} \right) \right)$$
and 
$$p_{q} = \frac{k}{\epsilon} \left( 1 + O \left( \frac{\ln \ln (k/\epsilon)}{\ln (k/\epsilon)} \right) \right).$$
Thus, if we use $K = c \cdot (k / \epsilon) \left \lfloor \log_{(k / \epsilon)} N \right \rfloor + 1$ in order to construct $\mathcal{M}_{(k/\epsilon),K}$ we will have
$$q+K = \frac{c \cdot k \left\lfloor \log_{(k/\epsilon)} N \right\rfloor}{\epsilon} \left( 1 + O \left( \frac{1}{\ln N} \right) \right).$$
Here we have assumed that $(k/\epsilon)+K$ is less than $N$.  Applying the Prime Number Theorem once more we have that
\begin{equation}
p_{q+K} = \frac{c \cdot k \left\lfloor \log_{(k/\epsilon)} N \right\rfloor \cdot \ln \left(\frac{k \cdot \ln N}{\epsilon} \right)}{\epsilon} \left( 1 + O\left( \frac{\ln \ln \left( \frac{k \ln N}{\epsilon} \right)}{\ln \left( \frac{k \ln N}{\epsilon} \right) } \right) \right).
\label{equ:pq+K}
\end{equation}
Utilizing Equation~\ref{eqn:Sum} now yields
\begin{equation}
m \leq \sum^{K-1}_{j=0} p_{q+j} = \frac{c^2 \cdot k^2 \left\lfloor \log_{(k/\epsilon)} N \right\rfloor^2 \cdot \ln \left(\frac{k \cdot \ln N}{\epsilon} \right)}{2 \epsilon^2} \left( 1 + O\left( \frac{\ln \ln \left( \frac{k \ln N}{\epsilon} \right)}{\ln \left( \frac{k \ln N}{\epsilon} \right) } \right) \right).
\label{eqn:NumRows}
\end{equation}
Hence, we have an asymptotic upper bound for the number of rows in $\mathcal{M}_{(k/\epsilon),K}$.  The next theorem, proven in Appendix~\ref{app:proof_thm:RowBound}, provides a concrete upper bound.

\begin{thm}
Suppose that $N, k, \epsilon^{-1} \in \mathbbm{N}-\{ 1 \}$ with $N > k \geq 2$.  Then, if we set $K = c \cdot (k / \epsilon) \left \lfloor \log_{s_1} N \right \rfloor + 1$ for some constant integer $c \in [2,\infty)$, there exists an $m \times N$ measurement matrix, $\mathcal{M}_{s_1,K}$, with a number of rows
$$m < \frac{3 (c + 1.89)^2 \cdot k^2 \left\lfloor \log_{(k/\epsilon)} N \right\rfloor^2}{4 \cdot \epsilon^2} \cdot \ln\left( \frac{(c + 1.89) \cdot k \left\lfloor \log_{(k/\epsilon)} N \right\rfloor}{\epsilon} \right).$$
Tighter upper bounds for the number of rows may be explicitly calculated using Equations~\ref{eqn:q+Kbound} -- \ref{eqn:UpperSum} below. 
\label{thm:RowBound}
\end{thm}

\noindent \textit{Proof:}  See Appendix~\ref{app:proof_thm:RowBound}.~~$\Box$ \\

Theorems~\ref{thm:Mest} and~\ref{thm:RowBound} collectively provide bounds for the number of rows a measurement matrix $\mathcal{M}_{s_1,K}$ may contain and still be able to estimate any entry of a vector $\vec{x} \in \mathbbm{C}^N$ to within a precision proportional to $\left\| \vec{x} - \vec{x}^{\rm opt}_{(k/\epsilon)} \right\|_1$.  These bounds are universal in that they pertain to measurement matrices which are guaranteed to provide accurate estimates for all entries of all vectors $\vec{x} \in \mathbbm{C}^N$.  In the next section we will prove the existence of a small number of $\mathcal{M}_{s_1,K}$ rows which are guaranteed to provide precise estimates for any sufficiently small set of vector entries.  We will also briefly consider a randomized matrix construction based on uniformly sampling rows of the deterministic $\mathcal{M}_{s_1,K}$ matrices considered above.  These results will ultimately motivate the development of sparse Fourier transforms with reduced sampling requirements.

\subsection{Randomized Row Sampling and Existence Results}
\label{sec:FixedSig}

In this section we will consider submatrices of the $m \times N$ measurement matrices, $\mathcal{M}_{s_1,K}$, discussed above.  More specifically, we will be discussing matrices formed by selecting a small number of rows from an $\mathcal{M}_{s_1,K}$ matrix as follows.  Let $\tilde{S} = \left\{ s_{j_1}, s_{j_2}, \dots, s_{j_l} \right\}$ be a subset of the $s_j$ values used to form $\mathcal{M}_{s_1,K}$ (see Equations~\ref{def:sj} --~\ref{eqn:Def_M}).  We will then define $\mathcal{M}_{\tilde{S}}$ to be the $\left( \tilde{m} = \sum^l_{\tilde{h}=1} s_{j_{\tilde{h}}} \right) \times N$ matrix,
\begin{equation}
\mathcal{M}_{\tilde{S}} = \left( \begin{array}{l} r_{j_1,0} \\ r_{j_1,1}\\ \vdots \\ r_{j_1,s_{j_1}-1} \\ \vdots \\ r_{j_l,s_{j_l}-1} \\
\end{array} \right),
\label{eqn:Def_MS}
\end{equation}
with each row defined as per Equation~\ref{eqn:Def_r}.  Finally, for $n \in [0,N) \cap \mathbbm{N}$, we define
$\mathcal{M}_{\tilde{S},n}$ to be the $l \times N$ matrix,
\begin{equation}
\mathcal{M}_{\tilde{S},n} = \left( \begin{array}{l} r_{j_1,n \textrm{~mod~} s_{j_1}} \\ 
r_{j_2,n \textrm{~mod~} s_{j_2}}\\ \vdots \\ r_{j_l,n \textrm{~mod~} s_{j_l}}\\ 
\end{array} \right),
\label{eqn:Def_MSn}
\end{equation}
along the lines of Equation~\ref{def:Msub}.  The following corollary of Theorem~\ref{thm:Mest} demonstrates the existence of small submatrices of $\mathcal{M}_{s_1,K}$ capable of providing accurate approximations to any given subset of a given vector $\vec{x} \in \mathbbm{C}^{N}$.

\begin{cor}
Let $k,N, \epsilon^{-1} \in \mathbbm{N}$, $S \subseteq [0,N) \cap \mathbbm{N}$, and $\vec{x} \in \mathbbm{C}^{N}$.  Set $K = c \cdot (k / \epsilon) \left \lfloor \log_{s_1} N \right \rfloor + 1$ for $s_1 \in \mathbbm{N}$ and a constant integer $c \in [4,\infty)$.  Form an $m \times N$ measurement matrix $\mathcal{M}_{s_1,K}$ as per Section~\ref{sec:makeM}.  Then, there exists a subset of $O\left(\log |S| \right)$ $s_j$ values for $\mathcal{M}_{s_1,K}$,
$$\tilde{S} = \left\{s_{j_1}, s_{j_2}, \dots, s_{j_{\left\lceil \log_{(c/2)} (|S|+1) \right\rceil}} \right\},$$
with the following property:  For all $n \in S$ we have
$$\min_{s_{j_h} \in \tilde{S}} \left| \left(\mathcal{M}_{\tilde{S},n} \vec{x} - x_n \cdot \vec{\mathbbm{1}}_{\left\lceil \log_{{(c/2)}} (|S|+1) \right\rceil} \right)_h \right| \leq \frac{ \epsilon \cdot \left\| \vec{x} - \vec{x}^{\rm opt}_{(k/\epsilon)} \right\|_1}{k}.$$
\label{cor:Exist}
\end{cor}

\noindent \textit{Proof:}\\

We proceed by induction on the size of $S \subseteq [0,N) \cap \mathbbm{N}$.  For the base case we assume $| S | = 1$ and apply Theorem~\ref{thm:Mest} with $n$ set to the single element of $S$.  We then define $\tilde{S}$ to be a singleton set containing any one of the $s_j$ rows of $\mathcal{M}_{s_1,K,n}$ which approximates $x_n$ to the guaranteed precision.  Now, suppose that the statement of Corollary~\ref{cor:Exist} holds for all subsets $S \subseteq [0,N) \cap \mathbbm{N}$ with $|S| \leq a \in \mathbbm{N}^{+}$.  Let $S' \subseteq [0,N) \cap \mathbbm{N}$ have $|S'| \leq \frac{a \cdot c}{2}$.  We will prove that the statement of Corollary~\ref{cor:Exist} holds for $S'$.

For each $n \in S'$ and $j \in [1,K] \cap \mathbbm{N}$ we will count a `failure' if $$\left| \left(\mathcal{M}_{s_1,K,n} \vec{x} \right)_j - x_n \right| > \frac{ \epsilon \cdot \left\| \vec{x} - \vec{x}^{\rm opt}_{(k/\epsilon)} \right\|_1}{k}.$$  Theorem~\ref{thm:Mest} tells us that there will be fewer than $(2/c) \cdot K$ `failures' for each element of $S'$, for a total of fewer than $\frac{2 \cdot |S'|}{c} \cdot K$ collective `failures' for all elements of $S'$.  Clearly, at least one of the $K$ $s_j$ values used to construct $\mathcal{M}_{s_1,K}$ must `fail' for fewer than $\frac{2 \cdot |S'|}{c}$ elements of $S'$.  Let $s'_j$ be the $s_j$ value which `fails' for the smallest number of elements of $S'$, and let $S'' \subset S'$ contain all the elements of $S'$ for which $s'_j$ `fails'.  We can see that $|S''| < \frac{2 \cdot |S'|}{c} \leq a$.  Our induction hypothesis applied to $S''$ together with the presence of $s'_j$ yields the desired result.~~$\Box$ \\  

Corollary~\ref{cor:Exist} demonstrates the existence of a small number of $s_j$ values which allow us to estimate every entry of a given vector.  However, it is apparently difficult to locate these $s_j$ values efficiently.  The following corollary circumvents this difficulty by showing that a small set of randomly selected $s_j$ values will still allow us to estimate all entries of any given vector with high probability.  Thus, in practice it suffices to select a random subset of the rows from a $\mathcal{M}_{s_1,K}$ matrix.

\begin{cor}
Let $k,N, \epsilon^{-1} \in \mathbbm{N}^{+}$, $\sigma \in [2/3,1)$, $S \subseteq [0,N) \cap \mathbbm{N}$, and $\vec{x} \in \mathbbm{C}^{N}$.  Set 
$K = c \cdot (k / \epsilon) \left \lfloor \log_{s_1} N \right \rfloor + 1$
for $s_1 \in \mathbbm{N}$ and a constant integer $c \in [14,\infty)$.  Form an $m \times N$ measurement matrix $\mathcal{M}_{s_1,K}$ as per Section~\ref{sec:makeM}.  Finally, form a multiset of the $s_j$ values for $\mathcal{M}_{s_1,K}$ by independently choosing 
\begin{equation}
l = \left\lceil 21 \cdot \ln \left( \frac{|S|}{1 - \sigma} \right) \right\rceil
\label{eqn:RandSj_lval}
\end{equation}
$s_j$ values uniformly at random with replacement.
Denote this multiset of $s_j$ values by
$$\tilde{S} = \left\{ s_{j_1}, s_{j_2}, \dots, s_{j_l} \right\}.$$
Then, with probability at least $\sigma$ the resulting random matrix, $\mathcal{M}_{\tilde{S}}$, will have the following property:  For all $n \in S$ more than $l / 2$ of the $s_{j_h} \in \tilde{S}$ (counted with multiplicity) will have
$$\left| \left(\mathcal{M}_{\tilde{S},n} \vec{x} - x_n \cdot \vec{\mathbbm{1}}_{l} \right)_h \right| \leq \frac{ \epsilon \cdot \left\| \vec{x} - \vec{x}^{\rm opt}_{(k/\epsilon)} \right\|_1}{k}.$$
\label{cor:RandSj}
\end{cor}

\noindent \textit{Proof:}  See Appendix~\ref{app:proof_cor:RandSj}.~~$\Box$ \\

Notice that Corollary~\ref{cor:RandSj} considers selecting a multiset of rows from a $\mathcal{M}_{s_1,K}$ measurement matrix.  In other words, some rows of the measurement matrix may be selected multiple times.  If this occurs in practice, one should consider any multiply selected rows to be chosen more than once for counting purposes only.  For example, during matrix multiplication a multiply selected row should be processed only once in order to avoid duplication of labor.  However, the results of these calculations should be considered multiple times for the purposes of estimation (e.g., in the median operations of Algorithm~\ref{alg:reconstruct1}).  

We will now consider these $m \times N$ matrices, $\mathcal{M}_{s_1,K}$, with respect to the discrete Fourier transform.  In particular, we will consider using $\mathcal{M}_{s_1,K}$ to estimate the Fourier transform of a periodic function along the lines of Theorem~\ref{thm:Mest}.  As we shall see, the special number theoretic nature of our matrix constructions will allow us to estimate Fourier coefficients of any periodic function by using a small number of function samples.

\subsection{The Fourier Case}
\label{sec:Fourier}

Suppose $f: [0, 2\pi] \rightarrow \mathbbm{C}$ is a complex valued function with $\hat{f} \in l^{1}$.  Let $P$ be the least common multiple of $\big\{N, s_1, \dots , s_K \big\}$ and form a set of samples from $f$, $\vec{A} \in \mathbbm{C}^P$, with
$$A_p = f\left( p \cdot \frac{2 \pi}{P} \right) \textrm{ for } p \in [0,P) \cap \mathbbm{N}.$$ 
Ultimately, we want to use $\mathcal{M}_{s_1,K} \vec{\hat{f}}$ in order to estimate the entries of the $N$-length vector $\vec{\hat{f}}$.  However, we must first calculate $\mathcal{M}_{s_1,K} \vec{\hat{f}}$.  In the remainder of this section we will discuss how to calculate $\mathcal{M}_{s_1,K} \vec{\hat{f}} \in \mathbbm{C}^m$ while using as few samples from $f$ as possible in the process.

To solve this problem we will use an extended version of our $m \times N$ matrix $\mathcal{M}_{s_1,K}$.  This extended matrix, $\mathcal{E}_{s_1,K}$, is the $m \times P$ matrix formed by extending each row $r_{j,h}$ of $\mathcal{M}_{s_1,K}$ as per Equation~\ref{eqn:Def_r} for all $p \in [0,P)$.  We now consider the product of $\mathcal{E}_{s_1,K}$ and the $P \times P$ discrete Fourier transform matrix, $\tilde{\Psi}$, defined by $\tilde{\Psi}_{\omega,p} = \frac{1}{P} \cdot e^{\frac{-2 \pi \mathbbm{i} \cdot \omega \cdot p}{P}}$.  For each row $r_{j,h}$ of $\mathcal{E}_{s_1,K}$ and column $p$ of $\tilde{\Psi}$ we have 
\begin{equation}
\left( \mathcal{E}_{s_1,K} \cdot \tilde{\Psi} \right)_{r_{j,h},p} ~~=~~ \frac{1}{P} \sum^{\frac{P}{s_j}-1}_{l=0} e^{\frac{-2 \pi \mathbbm{i} \cdot p \cdot \left( h + l \cdot s_j \right)}{P}} ~~=~~ \frac{e^{\frac{-2 \pi \mathbbm{i} \cdot p \cdot h}{P}}}{P} \sum^{\frac{P}{s_j}-1}_{l=0} e^{\frac{-2 \pi \mathbbm{i} \cdot p \cdot l}{P/s_{j}}} ~~=~~ \left\{ \begin{array}{ll} \frac{e^{\frac{-2 \pi \mathbbm{i} \cdot p \cdot h}{P}}}{s_{j}} & \textrm{if } p \equiv 0 \textrm{ mod } \frac{P}{s_{j}} \\ 0 & {\rm otherwise} \end{array} \right..
\label{eqn:FastMult}
\end{equation}

Thus, $\mathcal{E}_{s_1,K} \cdot \tilde{\Psi}$ is highly sparse.  In fact, we can see that each $r_{j,h}$ row contains only $s_{j}$ nonzero entries.  Better still, all the rows associated with a given $s_{j}$ have nonzero column entries in a pattern consistent with a small fast Fourier transform.  This aliasing phenomena results in a fast algorithm for computing $\mathcal{E}_{s_1,K} \cdot \tilde{\Psi} \cdot \vec{A}$ (see Algorithm~\ref{alg:fastMult}).  Lemma~\ref{lem:MultError} shows that $\mathcal{E}_{s_1,K} \tilde{\Psi} \vec{A}$ is a good approximation to $\mathcal{M}_{s_1,K} \vec{\hat{f}} \in \mathbbm{C}^m$ for all periodic functions whose Fourier transforms decay quickly enough.

\begin{algorithm}[tb]
\begin{algorithmic}[1]
\caption{$\proc{Fast Multiply}$} \label{alg:fastMult}
\STATE \textbf{Input: Function $f$, integers $k < K < N$, relatively prime $s_1, \dots, s_K$} 
\STATE \textbf{Output: $\mathcal{E}_{s_1,K} \cdot \tilde{\Psi} \cdot \vec{A}$}
\FOR {$j$ from $1$ to $K$}
	\STATE $\vec{A_{s_{j}}} \leftarrow f(0), f \left( \frac{2 \pi}{s_{j}} \right), \dots, f \left( \frac{2 \pi (s_{j} - 1)}{s_{j}} \right)$
	\STATE $\widehat{\vec{A_{s_{j}}}} \leftarrow$ \textbf{FFT} $\left[\vec{A_{s_{j}}} \right]$  
\ENDFOR 
\STATE Output $\left( \widehat{\vec{A_{s_{1}}}}, \widehat{\vec{A_{s_{2}}}}, \dots, \widehat{\vec{A_{s_{K}}}} \right)^{\textrm{T}}$
\end{algorithmic}
\end{algorithm}

\begin{lem}
Every entry of $\mathcal{E}_{s_1,K} \tilde{\Psi} \vec{A}$ approximates the associated entry of $\mathcal{M}_{s_1,K} \vec{\hat{f}}$ to within $\left\| ~\hat{f} - \bar{\hat{f}} ~ \right\|_1$ precision.
\label{lem:MultError}
\end{lem}

\noindent \textit{Proof:}\\  

Suppose that $N$ is odd (the case for $N$ even is analogous).  Then, for all $j \in [1,K] \cap \mathbbm{N}$ and $h \in [0,s_j) \cap \mathbbm{N}$, we have that
$$\left| \left( \mathcal{M}_{s_1,K} \vec{\hat{f}} - \mathcal{E}_{s_1,K} \tilde{\Psi} \vec{A} \right)_{r_{j,h}} \right| = \left| \sum_{l,~|h+l\cdot s_{j}| \leq \frac{N-1}{2}} \hat{f}_{h + l \cdot s_{j}} - \sum_{\omega \equiv h \textrm{ mod } s_{j}} \hat{f}(\omega) \right| = \left| \sum_{l,~|h+l\cdot s_{j}| \leq \frac{N-1}{2}} \hat{f} 
\left(h+l\cdot s_{j} \right) - \sum_{\omega \equiv h \textrm{ mod } s_{j}} \hat{f}(\omega) \right|.$$
Cancelling all Fourier coefficients for frequencies in $\left( - \big\lceil \frac{N}{2} \big\rceil, \big\lfloor \frac{N}{2} \big\rfloor \right] \cap \mathbbm{N}$ we get that
\begin{equation}
\left| \left( \mathcal{M}_{s_1,K} \vec{\hat{f}} - \mathcal{E}_{s_1,K} \tilde{\Psi} \vec{A} \right)_{r_{j,h}} \right| = \left| \sum_{\omega \equiv h \textrm{ mod } s_{j},~ |\omega| \geq \frac{N+1}{2}} \hat{f}(\omega) \right| \leq \sum_{|\omega| \geq \frac{N+1}{2}}  \left| \hat{f}(\omega) \right| = \left\|~ \hat{f} - \bar{\hat{f}} ~\right\|_1.~~\Box
\label{eqn:LemBound}
\end{equation}

By inspecting Equation~\ref{eqn:FastMult} it is not difficult to see that Algorithm~\ref{alg:fastMult} utilizes exactly $m - \left(K-1 \right)$ samples from $f$.
Considering this in combination with Theorem~\ref{thm:RowBound} in Section~\ref{sec:makeM} leads us to the conclusion that Algorithm~\ref{alg:fastMult} samples $f$ at $O\left( \frac{k^2 \cdot \left\lfloor \log_{(k/\epsilon)} N \right\rfloor^2 \cdot \ln \left(\frac{k \cdot \ln N}{\epsilon} \right)}{\epsilon^2} \right)$ distinct values.  Similarly, we can see that Algorithm~\ref{alg:fastMult} runs in time $O\left( \sum_{j=1}^{K} s_{j} \log s_{j} \right)$ if we calculate the FFTs using a chirp $z$-transform \cite{rabiner-schafer-rader}.  Thus, for well chosen $s_j$ values the runtime will be
\begin{align}
O\left( \sum_{j=1}^{K} s_{j} \log s_{j} \right) = O\left( \sum_{j=0}^{K-1} p_{q+j} \log p_{q+j} \right) & = O \left( p^2_{q+K} \right)~~\textrm{  (see~\cite{IDSFA})} \nonumber \\ & = O \left( \frac{k^2 \cdot \lfloor \log_{(k/\epsilon)} N \rfloor^2 \cdot \ln^2 \left(\frac{k \cdot \ln N}{\epsilon} \right)}{\epsilon^2} \right) 
\label{eqn:Alg1RUNTIME}
\end{align}
using Equation~\ref{equ:pq+K}.  We will now demonstrate how the specialized $m \times N$ matrices, $\mathcal{M}_{s_1,K}$, along with their extended $m \times P$ counterpart matrices, $\mathcal{E}_{s_1,K}$, considered throughout Sections~\ref{sec:makeM} and~\ref{sec:Fourier} can be utilized to construct accurate sparse Fourier transform methods.

\section{Fourier Reconstruction}
\label{sec:FourierRecov}

In this section we develop a sparse Fourier transform based on the measurement matrices considered in the previous section.  This sparse Fourier method is entirely dependent on the ability of our developed measurement matrices to accurately estimate any entry of a vector with which they have been multiplied (i.e., Theorem~\ref{thm:Mest}).  The idea behind the algorithm is simple.  We first quickly approximate the product of a Section~\ref{sec:makeM} measurement matrix with the Fourier transform of an input function using Algorithm~\ref{alg:fastMult}.  We then use the this product to accurately estimate all Fourier entries, keeping only the largest magnitude estimates for our final sparse Fourier approximation.  See Algorithm~\ref{alg:reconstruct1} for pseudo code.  Theorem~\ref{thm:Recon1} provides error, sampling, and runtime bounds for Algorithm~\ref{alg:reconstruct1}.

\begin{algorithm}[tb]
\begin{algorithmic}[1]
\caption{$\proc{Fourier Approximate 1}$} \label{alg:reconstruct1}
\STATE \textbf{Input: $k, N, \epsilon^{-1} \in \mathbbm{N}-\{ 1\}$}, Function $f$, Measurement matrix $\mathcal{M}_{s_1,K}$ with $K = 4 \cdot (k / \epsilon) \left \lfloor \log_{s_1} N \right \rfloor + 1$ (see Section~\ref{sec:makeM})
\STATE \textbf{Output: $\vec{x}_S$, an approximation to } $\vec{\hat{f}}^{\rm~ opt}_k$
\STATE Initialize $S \leftarrow \emptyset, ~\vec{x} \leftarrow \vec{0}_N$
\STATE $\mathcal{E}_{s_1,K} \tilde{\Psi} \vec{A} \leftarrow $ Algorithm~\ref{alg:fastMult}($f,~k,~K,~N,~s_j$ values for $\mathcal{M}_{s_1,K}$) 
\FOR {$\omega$ from $1 - \big\lceil \frac{N}{2} \big\rceil$ to $\big\lfloor \frac{N}{2} \big\rfloor$}
	\STATE $\mathbbm{Re}\left\{ x_{\omega} \right\} \leftarrow \rm{median~of~multiset} \left\{ \mathbbm{Re} \left\{ \left(\mathcal{E}_{s_1,K,\omega} \tilde{\Psi} \vec{A} \right)_j \right\}~\big|~1 \leq j \leq K \right\}$
	\STATE $\mathbbm{Im}\left\{ x_{\omega} \right\} \leftarrow \rm{median~of~multiset} \left\{ \mathbbm{Im} \left\{ \left(\mathcal{E}_{s_1,K,\omega} \tilde{\Psi} \vec{A} \right)_j \right\}~\big|~1 \leq j \leq K \right\}$
\ENDFOR
\STATE Sort $\vec{x}$ entries by magnitude so that $|x_{\omega_1}| \geq |x_{\omega_2}| \geq |x_{\omega_3}| \geq \dots$
\STATE $S \leftarrow \{ \omega_1, \omega_2, \dots, \omega_{2k} \}$
\STATE Output $\vec{x}_{S}$
\end{algorithmic}
\end{algorithm}

\begin{thm}
Suppose $f: [0,2\pi] \rightarrow \mathbbm{C}$ has $\hat{f} \in l^1$. Let $N, k, \epsilon^{-1} \in \mathbbm{N}-\{ 1 \}$ with $N > (k/\epsilon) \geq 2$.  Then, Algorithm~\ref{alg:reconstruct1} will output an $\vec{x}_{S} \in \mathbbm{C}^N$ satisfying
\begin{equation}
\left\| \vec{\hat{f}} - \vec{x}_{S} \right\|_2 \leq \left\| \vec{\hat{f}} - \vec{\hat{f}}^{\rm~ opt}_k \right\|_2 + \frac{22\epsilon \cdot \left\| \vec{\hat{f}} - \vec{\hat{f}}^{\rm~ opt}_{(k/\epsilon)} \right\|_1}{\sqrt{k}} + 22\sqrt{k} \cdot \left\| ~\hat{f} - \bar{\hat{f}} ~ \right\|_1.
\label{eqn:Recon1Error}
\end{equation} 
In the process $f$ will be evaluated at less than 
$$26.02 \cdot \frac{k^2 \left\lfloor \log_{(k/\epsilon)} N \right\rfloor^2}{\epsilon^2} \cdot \ln\left( \frac{5.89 \cdot k \left\lfloor \log_{(k/\epsilon)} N \right\rfloor}{\epsilon} \right)$$
points in $[0,2\pi]$. The runtime of lines 5 through 11 is $O\left( N \cdot (k / \epsilon) \log_{(k/\epsilon)} N \right)$.
\label{thm:Recon1}
\end{thm}

\noindent \textit{Proof:}\\

Fix $\omega \in \left( - \left\lceil \frac{N}{2} \right\rceil, \left\lfloor \frac{N}{2} \right\rfloor \right] \cap \mathbbm{Z}$ and let $\delta$ be set to  
$$\delta = \frac{ \epsilon \cdot \left\| \vec{\hat{f}} - \vec{\hat{f}}^{\rm ~opt}_{(k/\epsilon)} \right\|_1}{k} + \left\| ~\hat{f} - \bar{\hat{f}} ~ \right\|_1.$$
As a consequence of Theorem~\ref{thm:Mest} and Lemma~\ref{lem:MultError} we can see than more than half of the $K = 4 \cdot (k / \epsilon) \left \lfloor \log_{s_1} N \right \rfloor + 1$ entries of $\mathcal{E}_{s_1,K,\omega} \tilde{\Psi} \vec{A}$ produced in line 4 will satisfy $\left| \left(\mathcal{E}_{s_1,K,\omega} \tilde{\Psi} \vec{A} \right)_j - \hat{f}_\omega \right| \leq \delta$.  Therefore, the $x_{\omega}$ value produced by lines 6 and 7 will have 
\begin{equation}
\left| x_{\omega} - \hat{f}_\omega \right| \leq \sqrt{2} \cdot \delta.
\label{eqn:ApproxE}
\end{equation}

Since Equation~\ref{eqn:ApproxE} holds for all $\omega \in \left( - \left\lceil \frac{N}{2} \right\rceil, \left\lfloor \frac{N}{2} \right\rfloor \right] \cap \mathbbm{Z}$ we can begin to bound the approximation error by 
\begin{align}
\left\| \vec{\hat{f}} - \vec{x}_{S} \right\|_2 &~\leq~ \left\| \vec{\hat{f}} - \vec{\hat{f}}_{S} \right\|_2 + \left\| \vec{\hat{f}}_{S} - \vec{x}_{S} \right\|_2 ~\leq~ \left\| \vec{\hat{f}} - \vec{\hat{f}}_{S} \right\|_2 + 2\sqrt{k} \cdot \delta \nonumber \\
&~=~\sqrt{\left\| \vec{\hat{f}} - \vec{\hat{f}}^{\rm~ opt}_k \right\|^2_2 + \sum_{\omega \in S^{\rm opt}_k - S} \left| \hat{f}_{\omega} \right|^2 - \sum_{\tilde{\omega} \in S - S^{\rm opt}_k } \left| \hat{f}_{\tilde{\omega}} \right|^2} + 2\sqrt{k} \cdot \delta.
\label{eqn:ApproxEE}
\end{align}
In order to make additional progress on Equation~\ref{eqn:ApproxEE} we must first consider the possible magnitudes of $\vec{\hat{f}}$ entries at indices in $S - S^{\rm opt}_k$ and $S^{\rm opt}_k - S$.

Suppose $\omega \in S^{\rm opt}_k - S \neq \emptyset$ and let $\tilde{\omega} \in S - S^{\rm opt}_k$.  Line 10 will only have placed $\tilde{\omega} \in S$ instead of $\omega$ if $|x_{\tilde{\omega}}| \geq |x_{\omega}|$.  However, this can only happen if 
$$\left|\hat{f}_{\omega_{k}} \right| + \sqrt{2} \cdot \delta ~\geq~ \left|\hat{f}_{\tilde{\omega}} \right| + \sqrt{2} \cdot \delta ~\geq~ \left|\hat{f}_{\omega} \right| - \sqrt{2} \cdot \delta ~\geq~ \left|\hat{f}_{\omega_{k}} \right| - \sqrt{2} \cdot \delta.$$
In other words, all elements of $S - S^{\rm opt}_k$ and $S^{\rm opt}_k - S$ must index $\vec{\hat{f}}$ entries with roughly the same magnitude as the $k^{\rm th}$ largest magnitude entry of $\vec{\hat{f}}$ (up to a $\delta$ factor).  Furthermore, since $|S| = 2k$ we can see that $|S - S^{\rm opt}_k| \geq 2 \cdot |S^{\rm opt}_k - S|$.  We are now ready to give Equation~\ref{eqn:ApproxEE} further consideration.

If $S^{\rm opt}_k - S = \emptyset$ we are finished.  Otherwise, if $S^{\rm opt}_k - S \neq \emptyset$, we will have
$$\sum_{\tilde{\omega} \in S - S^{\rm opt}_k } \left| \hat{f}_{\tilde{\omega}} \right|^2 \geq 2 \cdot |S^{\rm opt}_k - S| \cdot \left( \left| \hat{f}_{\omega_{k}} \right| - 2 \sqrt{2} \cdot \delta \right)^2~=~{\bf A},$$
and
$${\bf B} = |S^{\rm opt}_k - S| \cdot \left( \left| \hat{f}_{\omega_{k}} \right| + 2 \sqrt{2} \cdot \delta \right)^2 \geq \sum_{\omega \in S^{\rm opt}_k - S} \left| \hat{f}_{\omega} \right|^2.$$
If ${\bf A} \geq {\bf B}$ then we are again finished.  If ${\bf A} < {\bf B}$ then 
$$\left| \hat{f}_{\omega_{k}} \right|^2 - 12 \sqrt{2} \delta \cdot \left| \hat{f}_{\omega_{k}} \right| + 8 \delta^2 < 0$$ 
which can only happen if $\left| \hat{f}_{\omega_{k}} \right| \in \left( (6 \sqrt{2} - 8) \cdot \delta, (6 \sqrt{2} + 8) \cdot \delta \right)$.  Therefore, in the worse case we can continue to bound Equation~\ref{eqn:ApproxEE} by
$$\left\| \vec{\hat{f}} - \vec{x}_{S} \right\|_2 \leq \sqrt{\left\| \vec{\hat{f}} - \vec{\hat{f}}^{\rm~ opt}_k \right\|^2_2 + k \cdot \left(8 \sqrt{2} + 8 \right)^2 \cdot \delta^2} ~+~ 2\sqrt{k} \cdot \delta ~\leq~ \left\| \vec{\hat{f}} - \vec{\hat{f}}^{\rm~ opt}_k \right\|_2 ~+~ 22\sqrt{k} \cdot \delta.$$
The error bound in Equation~\ref{eqn:Recon1Error} follows.

The upper bound on the number of point evaluations of $f$ follows directly from the application of Theorem~\ref{thm:RowBound} with $c = 4$.  Finding the largest $2k$ magnitude entries of $\vec{x}$ in lines 9 and 10 can be accomplished in $O(N \cdot \log k)$ time by using a binary search tree (see \cite{SortSearch}).  Therefore, the runtime of Algorithm~\ref{alg:reconstruct1} will be dominated by the median operations in lines 6 and 7.  Each of these medians can be accomplished in $O(K)$ time using a median-of-medians algorithm (e.g., \cite{MedianMedian}).  The stated $O(N \cdot K)$ runtime follows.~~$\Box$ \\

Note that the overall runtime behavior of Algorithm~\ref{alg:reconstruct1} will be dictated by both Equation~\ref{eqn:Alg1RUNTIME} and the runtime stated in Theorem~\ref{thm:Recon1}.  However, for most reasonable values of sublinear sparsity (i.e., whenever $k/\epsilon$ is $O(N / \log^3 N)$) the total runtime of Algorithm~\ref{alg:reconstruct1} will be $O\left( N \cdot (k / \epsilon) \log_{(k/\epsilon)} N \right)$.  One strategy for decreasing the runtime of Algorithm~\ref{alg:reconstruct1} is to decrease the number of measurement matrix rows, $K$, required to accurately estimate each Fourier coefficient.  Pursuing this strategy also has the additional benefit of reducing the number of function evaluations required for approximate Fourier reconstruction.  However, in exchange for these improvements we will have to sacrifice approximation guarantees for a small probability of outputting a relatively inaccurate answer.

Following the strategy above we will improve the performance of Algorithm~\ref{alg:reconstruct1} by modifying its input measurement matrix.  Instead of inputing a $\mathcal{M}_{s_1,K}$ measurement matrix as constructed in Section~\ref{sec:makeM} we will utilize a randomly constructed $\mathcal{M}_{\tilde{S}}$ measurement matrix as described in Section~\ref{sec:FixedSig}.  Corollary~\ref{cor:RandSj} ensures that such a randomly constructed $\mathcal{M}_{\tilde{S}}$ matrix will be likely to have all the properties of $\mathcal{M}_{s_1,K}$ matrices that Algorithm~\ref{alg:reconstruct1} needs.  Hence, with high probability we will achieve output from Algorithm~\ref{alg:reconstruct1} with the same approximation error bounds as derived for Theorem~\ref{thm:Recon1}.  Formalizing these ideas we obtain the following Corollary proved in Appendix~\ref{app:proof_cor:RandRecon1}.

\begin{cor}
Suppose $f: [0,2\pi] \rightarrow \mathbbm{C}$ has $\hat{f} \in l^1$. Let $\sigma \in [2/3,1)$ and $N, k, \epsilon^{-1} \in \mathbbm{N}-\{ 1 \}$ with $N > (k/\epsilon) \geq 2$.  Algorithm~\ref{alg:reconstruct1} may be executed using a matrix $\mathcal{M}_{\tilde{S}}$ from Section~\ref{sec:FixedSig} in place of the matrix $\mathcal{M}_{s_1,K}$ from Section~\ref{sec:makeM} to produce an output vector $\vec{x}_{S} \in \mathbbm{C}^N$ which will satisfy Equation~\ref{eqn:Recon1Error} with probability at least $\sigma$.  In the process $f$ will be evaluated at less than 
$$15.89 \cdot \left\lceil 21 \cdot \ln \left( \frac{N}{1 - \sigma} \right) \right\rceil \cdot \frac{k \left\lfloor \log_{(k/\epsilon)} N \right\rfloor }{\epsilon} \cdot \left( \ln \left( \frac{15.89 \cdot k \left\lfloor \log_{(k/\epsilon)} N \right\rfloor }{\epsilon} \right) + \ln \ln \left( \frac{15.89 \cdot k \left\lfloor \log_{(k/\epsilon)} N \right\rfloor }{\epsilon} \right) \right)$$
points in $[0,2\pi]$. The runtime of lines 5 through 11 will be $O\left( N \cdot \log \left( \frac{N}{1 - \sigma} \right) \right)$.
\label{cor:RandRecon1}
\end{cor}

\noindent \textit{Proof:}  See Appendix~\ref{app:proof_cor:RandRecon1}.~~$\Box$ \\

When executed with a random matrix $\mathcal{M}_{\tilde{S}}$ as input the overall runtime complexity of Algorithm~\ref{alg:reconstruct1} will be determined by both the runtime stated in Corollary~\ref{cor:RandRecon1} and the runtime of Algorithm~\ref{alg:fastMult}.  Suppose $\tilde{S}$ is a subset of $O\left( \log \left( \frac{N}{1 - \sigma} \right) \right)$ $s_j$ values defined as per Equations~\ref{eqn:Primes} --~\ref{eqn:Def_k}.  Then, Algorithm~\ref{alg:fastMult} will have a runtime complexity of
\begin{align}
O\left( \sum_{s_j \in \tilde{S}} s_{j} \cdot \log s_{j} \right) & = O\left( p_{q+K} \cdot \log p_{q+K} \cdot \log \left( \frac{N}{1 - \sigma} \right) \right) ~~\textrm{  (see~Equation~\ref{equ:pq+K})} \nonumber \\ & = 
O\left( \frac{k \cdot \log_{(k/\epsilon)} N}{\epsilon} \cdot \log^2 \left( \frac{k \cdot \log N }{\epsilon} \right) \cdot \log \left( \frac{N}{1 - \sigma} \right) \right).
\label{eqn:Alg1RandRUNTIME}
\end{align}
Thus, Algorithm~\ref{alg:reconstruct1} executed with a random input matrix from Section~\ref{sec:FixedSig} will have a total runtime complexity of $O\left( N \cdot \log \left( \frac{N}{1 - \sigma} \right) \right)$ whenever $(k/\epsilon)$ is $O(N / \log^3 N)$.  If we now set the desired success probability, $\sigma$, to be $1 - 1 / N^{O(1)}$ we obtain an overall $O(N \cdot \log N)$ computational complexity for Algorithm~\ref{alg:reconstruct1}.  This matches the runtime behavior of a standard fast Fourier transform while requiring asymptotically fewer function evaluations.

In the next section we will discuss methods for further decreasing the runtime requirements of Algorithm~\ref{alg:reconstruct1} while maintaining its approximation guarantees (i.e., the error bound in Equation~\ref{eqn:Recon1Error}).  
As a result we will develop sublinear-time Fourier algorithms that have both universal recovery guarantees and uniformly bounded runtime requirements.

\section{Decreasing the Runtime Complexity}
\label{sec:FourierRecov2}

Let $\mathcal{A},~\mathcal{B}$ be $m \times N$ and $\tilde{m} \times N$ complex valued matrixes, respectively.  Then, their row tensor product, $\mathcal{A} \circledast \mathcal{B}$, is defined to be the $\left( m \cdot \tilde{m} \right) \times N$ complex valued matrix created by performing component-wise multiplication of all rows of $\mathcal{A}$ with all rows of $\mathcal{B}$.  More specifically,
$$\left(\mathcal{A} \circledast \mathcal{B}\right)_{i,j} = \mathcal{A}_{i \textrm{~mod~m},j} \cdot \mathcal{B}_{\frac{i - i \textrm{~mod~m}}{m},j}.$$
In this section we will use the row tensor product of two types of specially constructed measurement matrices in order to improve the runtime complexity of Algorithm~\ref{alg:reconstruct1}.  One of these matrix types will be the $\mathcal{M}_{s_1,K}$ measurement matrices developed in Section~\ref{sec:makeM}.
The other type of matrix is described in the next two paragraphs.

Suppose that an $m \times N$ measurement matrix, $\mathcal{M}_{s_1,K}$, is given.  Furthermore, suppose that $s_1, \dots, s_K \in \mathbbm{N}$ are such that there exist $\lambda$ integers, $t_1 < \cdots < t_\lambda < s_1$, with $$\prod^\lambda_{i=1} t_i \geq \frac{N}{s_1}$$
that also have the property that the set
$$\left\{ t_1, \dots, t_\lambda, s_1, \dots, s_K \right\}$$
is pairwise relatively prime.
Note that such $t_i$ values can indeed be found if all the given $s_j$ values are prime numbers and $s_1 \geq \log_2 N \cdot \left( \ln \log_2 N + \ln \ln \log_2 N \right) \geq p_{\lfloor \log_2 N \rfloor}$ for $N \geq 64$ (see \cite{PrimeBound}).  We will now demonstrate how to use such $t_i$ values to create an $\tilde{m} \times N$ matrix, $\mathcal{N}_{\lambda,s_1}$, along the lines of Section~\ref{sec:makeM}.

Create a row, $\tilde{r}_{i,h}$, in $\mathcal{N}_{\lambda,s_1}$ for each possible residue of each $t_i$ integer (i.e., $\tilde{r}_{i,h}$ has $i \in [1,\lambda] \cap \mathbbm{N}$ and $h \in [0,t_i) \cap \mathbbm{N}$).  The $n^{\rm th}$ entry of each $\tilde{r}_{i,h}$ row, $n \in [0,N) \cap \mathbbm{N}$, will be
\begin{align}
(\tilde{r}_{i,h})_{n} & = \delta \left( (n-h) \textrm{ mod } t_i \right) \nonumber \\ & = \left\{ \begin{array}{ll} 1 & \textrm{if } n \equiv h \textrm{ mod } t_i \\ 0 & {\rm otherwise} \end{array} \right. .
\label{eqn:Def_r_tilde}
\end{align}
We then define
\begin{equation}
\mathcal{N}_{\lambda,s_1} = \left( \begin{array}{l} \vec{\mathbbm{1}}_N \\ \tilde{r}_{1,0} \\ \vdots \\ \tilde{r}_{1,t_{1}-1} \\ \vdots \\ \tilde{r}_{\lambda,t_{\lambda}-1} \\
\end{array} \right).
\label{eqn:Def_N}
\end{equation}
The result is an $\left( \tilde{m} = 1 + \sum^\lambda_{i = 1} t_i \right) \times N$ matrix with binary entries.  The following Lemma, proven in Appendix~\ref{app:proof_lem:RowBound}, upper bounds the smallest possible number of rows in any such $\mathcal{N}_{\lambda,s_1}$ matrix.

\begin{lem}
Suppose that $N, s_1, \dots, s_K \in \mathbbm{N}$ with 
$$\frac{N}{3} \geq s_1 > \left\lceil 3 \cdot \frac{\ln \left( N / s_1 \right)}{\ln \ln \left( N / s_1 \right)} \right\rceil \cdot \left( \ln \left\lceil 3 \cdot \frac{\ln \left( N / s_1 \right)}{\ln \ln \left( N / s_1 \right)} \right\rceil + \ln \ln \left\lceil 3 \cdot \frac{\ln \left( N / s_1 \right)}{\ln \ln \left( N / s_1 \right)} \right\rceil \right),$$ and $s_1, \dots, s_K$ containing no prime factors less than $s_1$.  Then, there exists a valid $\tilde{m} \times N$ measurement matrix, $\mathcal{N}_{\lambda,s_1}$, with a number of rows
$$\tilde{m} < \frac{3}{4}\left(\left\lceil 3 \cdot \frac{\ln \left( N / s_1 \right)}{\ln \ln \left( N / s_1 \right)} \right\rceil + 1\right)^2 \cdot \ln\left(\left\lceil 3 \cdot \frac{\ln \left( N / s_1 \right)}{\ln \ln \left( N / s_1 \right)} \right\rceil+1 \right) +1.$$
The corresponding value of $\lambda$ is $\lceil 3 \cdot \ln(N/s_1) / \ln \ln(N/s_1) \rceil$.
\label{lem:RowBound}
\end{lem}
\noindent \textit{Proof:}  See Appendix~\ref{app:proof_lem:RowBound}.~~$\Box$ \\

The $\left(m \cdot \tilde{m} \right) \times N$ row tensor product matrix, $\mathcal{R}_{\lambda,K} = \mathcal{M}_{s_1,K} \circledast \mathcal{N}_{\lambda,s_1}$, has several useful properties.  First, the fact that the first row of $\mathcal{N}_{\lambda,s_1}$ is the all-ones vector means that $\mathcal{R}_{\lambda,K}$ will contain a copy of every row of $\mathcal{M}_{s_1,K}$.  Second, all $\mathcal{R}_{\lambda,K}$ rows that are not copies of $\mathcal{M}_{s_1,K}$ rows will have the form $\bar{r}_{i,j,h} = r_{j,h~{\rm mod}~s_j} \circledast \tilde{r}_{i,h~{\rm mod}~t_i}$ for some $i \in [1,\lambda] \cap \mathbbm{N},~j \in [1,K] \cap \mathbbm{N}$, and $h \in [0,t_i \cdot s_j) \cap \mathbbm{N}$.  That is, the Chinese Remainder Theorem tells us that each such $\mathcal{R}_{\lambda,K}$ row will have its $n^{\rm th}$ entry given by 
\begin{align}
(\bar{r}_{i,j,h})_{n} = \delta \left( (n-h) \textrm{ mod } t_i \cdot s_j \right) = \left\{ \begin{array}{ll} 1 & \textrm{if } n \equiv h \textrm{ mod } t_i \cdot s_j \\ 0 & {\rm otherwise} \end{array} \right. .
\label{eqn:Def_r_bar}
\end{align}
The end result is that $\mathcal{R}_{\lambda,K}$ maintains a rigid number theoretic structure.  The following Lemma summarizes the most important properties of $\mathcal{R}_{\lambda,K} = \mathcal{M}_{s_1,K} \circledast \mathcal{N}_{\lambda,s_1}$.

\begin{lem}
Let $k, \epsilon^{-1}, s_1, \lambda,n \in [2,N) \cap \mathbbm{N}$, $\vec{x} \in \mathbbm{C}^N$, and $K = 4 \cdot (k / \epsilon) \left \lfloor \log_{s_1} N \right \rfloor + 1$.  Then, more than $\frac{K}{2}$ of the $K$ entries of $\mathcal{M}_{s_1,K,n} \cdot \vec{x}$ will estimate $x_{n}$ to within $\bar{\delta} = \frac{\epsilon \cdot \left\| \vec{x} - \vec{x}^{\rm opt}_{(k/\epsilon)} \right\|_1}{k}$ precision.  Furthermore, if $r_{j',n \textrm{~mod~} s_{j'}} \in \{ 0,~1 \}^N$ is a row of $\mathcal{M}_{s_1,K,n}$ associated with one of these more than $\frac{K}{2}$ entries then it will have all of the following properties: 
\begin{enumerate}
\item $ \left| r_{j',n \textrm{~mod~} s_{j'}} \cdot \vec{x} - x_{n} \right| \leq \bar{\delta}$,
\item $\left| \left( r_{j',n \textrm{~mod~} s_{j'}} \circledast \tilde{r}_{i, n \textrm{~mod~} t_{i}} \right) \cdot \vec{x} - x_{n} \right| = \left| \bar{r}_{i,j',n \textrm{~mod~} t_i \cdot s_{j'}} \cdot \vec{x} - x_{n} \right| \leq \bar{\delta}$ for all $i \in [1,\lambda] \cap \mathbbm{N}$, and
\item $\left| \left( r_{j',n \textrm{~mod~} s_{j'}} \circledast \tilde{r}_{i, h} \right) \cdot \vec{x} ~\right| = \left| \bar{r}_{i,j',\bar{h} \neq n \textrm{~mod~} t_i \cdot s_{j'}} \cdot \vec{x} ~\right| \leq \bar{\delta}$ for all $i \in [1,\lambda] \cap \mathbbm{N}$ and $h \in [0,t_i) \cap \left( \mathbbm{N} - \{ n \textrm{~mod~} t_i \} \right)$.
\end{enumerate}
\label{lem:Rprops}
\end{lem}
\noindent \textit{Proof:}  See Appendix~\ref{app:proof_lem:Rprops}.~~$\Box$\\

Suppose $f: [0, 2\pi] \rightarrow \mathbbm{C}$ is a complex valued function with $\hat{f} \in l^{1}$.  It is not difficult to see that $\mathcal{R}_{\lambda,K} \vec{\hat{f}}$ can be approximated using Algorithm~\ref{alg:fastMult} from Section~\ref{sec:Fourier} since $\mathcal{R}_{\lambda,K}$ maintains the required number theoretic structure.  We will simply perform FFTs on arrays of function samples with sizes given by all possible $t_i \cdot s_j$ value products.  The total number of function samples taken will be at most $m \cdot \tilde{m} - ( \lambda \cdot K + K - 1)$.  For $s_j$ and $t_i$ values chosen as per Theorem~\ref{thm:RowBound} and Lemma~\ref{lem:RowBound}, respectively, the runtime required by Algorithm~\ref{alg:fastMult} to approximate $\mathcal{R}_{\lambda,K}\vec{\hat{f}}$ will be
\begin{align}
O\left( \sum^{\lambda}_{i=1} \sum_{j=1}^{K} t_i \cdot s_{j} \log s_{j} \right) = O\left( \sum^{\lambda}_{i=1} \sum_{j=0}^{K-1} p_i \cdot p_{q+j} \log p_{q+j} \right) & = O \left( \frac{p^2_{q+K} \cdot p^2_{\lambda}}{\ln p_{\lambda}} \right)~~\textrm{  (see~\cite{IDSFA})} \nonumber \\ & = O \left( k^2 \cdot \frac{\ln^2 N \cdot \ln^2 \left(\frac{k \cdot \ln N}{\epsilon} \right) \cdot \ln^2 \left( \frac{\epsilon \cdot N}{k} \right)}{\epsilon^2 \cdot \ln^2 \left( \frac{k}{\epsilon} \right) \cdot \ln \ln \left( \frac{\epsilon \cdot N}{k}\right)}  \right). 
\label{eqn:Alg1RUNTIME2}
\end{align}
The last equality follows from Equation~\ref{equ:pq+K} and the Prime Number Theorem.  Finally, it is not difficult to see that the precision guarantees of  Lemma~\ref{lem:MultError} will still hold for an Algorithm~\ref{alg:fastMult} approximation to $\mathcal{R}_{\lambda,K} \vec{\hat{f}}$.  

Perhaps most importantly, the number theoretic structure of $\mathcal{R}_{\lambda,K}$ also allows us to use methods analogous to those outlined in Sections 1.1 and 5 of \cite{JournalDSFT} to quickly identify frequencies with large magnitude Fourier coefficients in $\hat{f}$.  Suppose that $\left| \hat{f}_{\omega} \right|$ is large relative to $\left\| \hat{f} \right\|_1$ (e.g., more than one tenth as large).  In this case Lemma~\ref{lem:Rprops} above tells us that $\hat{f}_{\omega}$ will also have a magnitude nearly as large as that of most entries of $\mathcal{M}_{s_1,K,\omega} \vec{\hat{f}}$.  Let $r_{j,\omega \textrm{~mod~} s_j}$ be the row of $\mathcal{M}_{s_1,K,\omega}$ associated with one of these $\mathcal{M}_{s_1,K,\omega} \vec{\hat{f}}$ entries dominated by $\hat{f}_{\omega}$.  By its construction we know that $\mathcal{R}_{\lambda,K}$ will not only contain $r_{j,\omega \textrm{~mod~} s_j}$, but also the related rows $\bar{r}_{1,j,\omega \textrm{~mod~} t_1 \cdot s_j},~\dots,~\bar{r}_{\lambda,j,\omega \textrm{~mod~} t_{\lambda} \cdot s_j}$.  Furthermore, all $\lambda + 1$ entries of $\mathcal{R}_{\lambda,K,\omega} \vec{\hat{f}}$ associated with these rows will also be dominated by $\hat{f}_{\omega}$ (see Lemma~\ref{lem:Rprops}).  On the other hand, for each $i \in [1,\lambda] \cap \mathbbm{N}$ the $\left(\mathcal{R}_{\lambda,K,\omega} \vec{\hat{f}}\right)_{\bar{r}_{i,j,h \neq \omega \textrm{~mod~} t_{i} \cdot s_j}}$ entries will all be significantly smaller than $\hat{f}_{\omega}$ in magnitude.  Hence, by comparing the relative magnitudes of the entries in $\left(r_{j,\omega \textrm{~mod~} s_j} \circledast \mathcal{N}_{\lambda,s_1}\right) \vec{\hat{f}}$ we can discern $\omega \textrm{~mod~} s_j, \omega \textrm{~mod~} t_1 \cdot s_j, \dots, \omega \textrm{~mod~} t_{\lambda} \cdot s_j$.  The end result is that $\omega$ can be recovered by inspecting $\mathcal{R}_{\lambda,K,\omega} \vec{\hat{f}}$.  See \cite{JournalDSFT} for a detailed discussion of a similar recovery procedure.  Utilizing these ideas we obtain Algorithm~\ref{alg:reconstruct2}.

\begin{algorithm}[tb]
\begin{algorithmic}[1]
\caption{$\proc{Fourier Approximate 2}$} \label{alg:reconstruct2}
\STATE \textbf{Input: $k, N, \epsilon^{-1} \in \mathbbm{N}-\{ 1\}$}, Function $f$, An $\left( m \cdot \tilde{m} \right) \times N$ measurement matrix $\mathcal{R}_{\lambda,K}$ with $K = 4 \cdot (k / \epsilon) \left \lfloor \log_{s_1} N \right \rfloor + 1$
\STATE \textbf{Output: $\vec{x}_S$, an approximation to } $\vec{\hat{f}}^{\rm~ opt}_k$
\STATE Initialize $S \leftarrow \emptyset, ~\vec{x} \leftarrow \vec{0}_N$
\STATE $\mathcal{G}_{\lambda,K} \tilde{\Psi} \vec{A} \leftarrow $ Algorithm~\ref{alg:fastMult}($f,~k,~K,~N,~s_j \textrm{ and } t_i$ values for $\mathcal{R}_{\lambda,K}$) 
\STATE $\mathcal{E}_{s_1,K} \tilde{\Psi} \vec{A} \leftarrow$ The $m$ entries of $\mathcal{G}_{\lambda,K} \tilde{\Psi} \vec{A}$ that approximate $\mathcal{M}_{s_1,K} \vec{\hat{f}}$ 
\begin{center}
{\sc Identification of Frequencies with Large Fourier Coefficients}
\end{center}
\FOR {$j$ from $1$ to $K$}
	\FOR {$h$ from $0$ to $s_{j}-1$}
		\FOR {$i$ from $1$ to $\lambda$}
			\STATE $b_{\rm min} \leftarrow \arg \min_{b \in [0,t_{i})} \left| \left(\mathcal{E}_{s_1,K} \tilde{\Psi} \vec{A}\right)_{r_{j,h}} - \left( \mathcal{G}_{\lambda,K} \tilde{\Psi} \vec{A} \right)_{\bar{r}_{i,j,h + b\cdot s_j}} \right|$
			\STATE $a_{j,h,i} \leftarrow \left( h + b_{\rm min} \cdot s_j \right)$ mod $t_i$
		\ENDFOR
		\STATE Reconstruct $\omega_{j,h}$ using that $\omega_{j,h} \equiv h \textrm{ mod } s_j, ~\omega_{j,h} \equiv a_{j,h,1} \textrm{ mod } t_1, ~\dots, ~\omega_{j,h} \equiv a_{j,h,\lambda} \textrm{ mod } t_\lambda$ 
	\ENDFOR
\ENDFOR \\
\begin{center}
{\sc Fourier Coefficient Estimation}
\end{center}
\FOR {\textbf{each} $\omega_{j,h}$ value reconstructed $> \frac{K}{2}$ times}
\STATE $\mathbbm{Re}\left\{ x_{\omega_{j,h}} \right\} \leftarrow \rm{median~of~multiset} \left\{ \mathbbm{Re} \left\{ \left(\mathcal{G}_{\lambda,K,\omega_{j,h}} \tilde{\Psi} \vec{A} \right)_j \right\}~\big|~1 \leq j \leq K \cdot \left( \lambda + 1 \right) \right\}$
	\STATE $\mathbbm{Im}\left\{ x_{\omega_{j,h}} \right\} \leftarrow \rm{median~of~multiset} \left\{ \mathbbm{Im} \left\{ \left(\mathcal{G}_{\lambda,K,\omega_{j,h}} \tilde{\Psi} \vec{A} \right)_j \right\}~\big|~1 \leq j \leq K \cdot \left( \lambda + 1 \right) \right\}$
\ENDFOR
\STATE Sort nonzero $\vec{x}$ entries by magnitude so that $|x_{\omega_1}| \geq |x_{\omega_2}| \geq |x_{\omega_3}| \geq \dots$
\STATE $S \leftarrow \{ \omega_1, \omega_2, \dots, \omega_{2k} \}$
\STATE Output $\vec{x}_{S}$
\end{algorithmic}
\end{algorithm}

Note that Algorithms~\ref{alg:reconstruct1} and~\ref{alg:reconstruct2} are quite similar.  The only significant difference between them is that Algorithm~\ref{alg:reconstruct1} estimates Fourier coefficients for all frequencies in the bandwidth specified by $N$ whereas Algorithm~\ref{alg:reconstruct2} restricts itself to estimating the Fourier coefficients for only a small number of frequencies it identifies as significant.  Given these similarities it should not be surprising that demonstrating the correctness of Algorithm~\ref{alg:reconstruct2} depends primarily on showing that it can correctly identify all frequencies with coefficients that are sufficiently large in magnitude.  This is established in Lemma~\ref{lem:FreqIdent} below.

\begin{lem}
Suppose that $\omega \in \left( - \left\lceil \frac{N}{2} \right\rceil, \left\lfloor \frac{N}{2} \right\rfloor \right] \cap \mathbbm{Z}$ is such that 
$$\left| \hat{f}_{\omega} \right| ~>~ 4 \cdot \left( \frac{ \epsilon \cdot \left\| \vec{\hat{f}} - \vec{\hat{f}}^{\rm ~opt}_{(k/\epsilon)} \right\|_1}{k} + \left\| ~\hat{f} - \bar{\hat{f}} ~ \right\|_1 \right).$$ 
Then, lines 6 through 14 of Algorithm~\ref{alg:reconstruct2} will reconstruct $\omega$ more than $\frac{K}{2}$ times.
\label{lem:FreqIdent}
\end{lem}
\noindent \textit{Proof:}\\

Suppose that $\omega \in \left( - \left\lceil \frac{N}{2} \right\rceil, \left\lfloor \frac{N}{2} \right\rfloor \right] \cap \mathbbm{Z}$ has $\left| \hat{f}_{\omega} \right| > 4 \delta$ where 
$$\delta = \frac{ \epsilon \cdot \left\| \vec{\hat{f}} - \vec{\hat{f}}^{\rm ~opt}_{(k/\epsilon)} \right\|_1}{k} + \left\| ~\hat{f} - \bar{\hat{f}} ~ \right\|_1.$$
Lemma~\ref{lem:Rprops} and Lemma~\ref{lem:MultError} guarantee that $\left| \left( \mathcal{E}_{s_1,K,\omega} \tilde{\Psi} \vec{A}\right)_{j} - \hat{f}_{\omega} \right| \leq \delta$ for more than $\frac{K}{2}$ entry indexes $j$.  Furthermore, if $j' \in [1,K] \cap \mathbbm{N}$ is one of these more than $\frac{K}{2}$ indexes, then property $(2)$ of Lemma~\ref{lem:Rprops} together with the preceding discussion of Lemma~\ref{lem:MultError} also ensures that $\left| \left( \mathcal{G}_{\lambda,K} \tilde{\Psi} \vec{A} \right)_{\bar{r}_{i,j',\omega ~{\rm mod}~ t_i \cdot s_{j'}}} - \hat{f}_{\omega} \right| \leq \delta$ for all $i \in [1,\lambda] \cap \mathbbm{N}$.  Fix $i \in [1,\lambda] \cap \mathbbm{N}$.  Thus, if $b \in [0,t_i) \cap \mathbbm{N}$ in line 9 of Algorithm~\ref{alg:reconstruct2} satisfies 
\begin{equation}
\omega \equiv \left( \left( \omega ~{\rm mod}~ s_{j'} \right) + b \cdot s_{j'} \right) ~{\rm mod}~ t_i \cdot s_{j'}
\label{eqn:bminSat}
\end{equation} 
we can see that 
$$\left| \left(\mathcal{E}_{s_1,K} \tilde{\Psi} \vec{A}\right)_{r_{j',\omega ~{\rm mod}~ s_{j'}}} - \left( \mathcal{G}_{\lambda,K} \tilde{\Psi} \vec{A} \right)_{\bar{r}_{i,j',\left( \omega ~{\rm mod}~ s_{j'} \right) + b\cdot s_{j'}}} \right| = \left| \left(\mathcal{E}_{s_1,K} \tilde{\Psi} \vec{A}\right)_{r_{j',\omega ~{\rm mod}~ s_{j'}}} - \hat{f}_{\omega} + \hat{f}_{\omega} - \left( \mathcal{G}_{\lambda,K} \tilde{\Psi} \vec{A} \right)_{\bar{r}_{i,j',\omega ~{\rm mod}~ t_i \cdot s_{j'}}} \right| \leq 2\delta.$$
Otherwise, if $b \in [0,t_i) \cap \mathbbm{N}$ does not satisfy Equation~\ref{eqn:bminSat}, property $(3)$ from Lemma~\ref{lem:Rprops} in combination with Lemma~\ref{lem:MultError} ensures that 
\begin{align}
2\delta &< \left| \hat{f}_{\omega} \right| - \left| \left(\mathcal{E}_{s_1,K} \tilde{\Psi} \vec{A}\right)_{r_{j',\omega ~{\rm mod}~ s_{j'}}} - \hat{f}_{\omega} \right| - \left| \left( \mathcal{G}_{\lambda,K} \tilde{\Psi} \vec{A} \right)_{\bar{r}_{i,j',\left( \omega ~{\rm mod}~ s_{j'} \right) + b\cdot s_{j'}}} \right| \nonumber \\
&\leq \left| \left(\mathcal{E}_{s_1,K} \tilde{\Psi} \vec{A}\right)_{r_{j',\omega ~{\rm mod}~ s_{j'}}} - \hat{f}_{\omega} + \hat{f}_{\omega} - \left( \mathcal{G}_{\lambda,K} \tilde{\Psi} \vec{A} \right)_{\bar{r}_{i,j',\left( \omega ~{\rm mod}~ s_{j'} \right) + b\cdot s_{j'}}} \right| = \left| \left(\mathcal{E}_{s_1,K} \tilde{\Psi} \vec{A}\right)_{r_{j',\omega ~{\rm mod}~ s_{j'}}} - \left( \mathcal{G}_{\lambda,K} \tilde{\Psi} \vec{A} \right)_{\bar{r}_{i,j',\left( \omega ~{\rm mod}~ s_{j'} \right) + b\cdot s_{j'}}} \right| \nonumber .
\end{align}
Therefore, the $b = b_{\rm min}$ identified in line 9 of Algorithm~\ref{alg:reconstruct2} will be guaranteed to satisfy Equation~\ref{eqn:bminSat} for all $i \in [1,\lambda] \cap \mathbbm{N}$.

Once we have identified $\omega ~{\rm mod}~ t_i \cdot s_{j'}$ in this fashion we can find $\omega ~{\rm mod}~ t_i$ in line 10 of Algorithm~\ref{alg:reconstruct2} by computing $\left( \omega ~{\rm mod}~ t_i \cdot s_{j'} \right) {\rm mod}~ t_i$.  Finally, by construction, the set $\{ t_1,~ \dots,~ t_\lambda,~ s_{j'} \}$ both has a collective product larger than $N$, and is pairwise relatively prime.  Therefore, the Chinese Remainder Theorem guarantees that line 12 of Algorithm~\ref{alg:reconstruct2} will indeed correctly reconstruct $\omega$ when $j = j'$ and $h = \omega ~{\rm mod}~ s_{j'}$.~~$\Box$ \\

With Lemma~\ref{lem:FreqIdent} in hand we are now prepared to prove that Algorithm~\ref{alg:reconstruct2} can indeed recover near-optimal sparse Fourier representations in sublinear-time.  We begin by using Lemma~\ref{lem:FreqIdent} to show that all sufficiently energetic frequencies are guaranteed to be identified.  Hence, the only way Algorithm~\ref{alg:reconstruct2} will not include an optimal Fourier representation frequency in its output is if the frequency is either $(i)$ insufficiently energetic to be identified, or $(ii)$ gets identified, but is then mistakenly estimated to have a smaller magnitude Fourier coefficient than many other somewhat energetic frequencies.  In the case of $(i)$ it is forgivable to exclude the frequency given that it must have a Fourier coefficient with a relatively small magnitude.  In the case of $(ii)$ we make up for the exclusion of a truly energetic frequency term by including many other less significant, but still fairly energetic, frequency terms in its place.  Carefully combining these ideas leads us to the error, sampling, and runtime bounds for Algorithm~\ref{alg:reconstruct2} stated in Theorem~\ref{thm:Recon2} below.

\begin{thm}
Suppose $f: [0,2\pi] \rightarrow \mathbbm{C}$ has $\hat{f} \in l^1$. Let $N, k, \epsilon^{-1} \in \mathbbm{N}-\{ 1 \}$ with $N > (k/\epsilon) \geq 2$.  Then, Algorithm~\ref{alg:reconstruct2} will output an $\vec{x}_{S} \in \mathbbm{C}^N$ satisfying
\begin{equation}
\left\| \vec{\hat{f}} - \vec{x}_{S} \right\|_2 \leq \left\| \vec{\hat{f}} - \vec{\hat{f}}^{\rm~ opt}_k \right\|_2 + \frac{22\epsilon \cdot \left\| \vec{\hat{f}} - \vec{\hat{f}}^{\rm~ opt}_{(k/\epsilon)} \right\|_1}{\sqrt{k}} + 22\sqrt{k} \cdot \left\| ~\hat{f} - \bar{\hat{f}} ~ \right\|_1.
\label{eqn:Recon2Error}
\end{equation} 
Under the conditions of Lemma~\ref{lem:RowBound}, $f$ will be evaluated at less than 
$$19.52 \cdot \frac{k^2 \left\lfloor \log_{(k/\epsilon)} N \right\rfloor^2}{\epsilon^2} \cdot \ln\left( \frac{5.89 \cdot k \left\lfloor \log_{(k/\epsilon)} N \right\rfloor}{\epsilon} \right) \cdot \left\{
\left(\left\lceil 3 \cdot \frac{\ln \left( \epsilon N/k \right)}{\ln \ln \left( \epsilon N/k \right)} \right\rceil + 1\right)^2 \cdot \ln\left(\left\lceil 3 \cdot \frac{\ln \left( \epsilon N/k \right)}{\ln \ln \left( \epsilon N/k \right)} \right\rceil+1 \right) + \frac{4}{3}\right\}$$
points in $[0,2\pi]$.  The runtime of lines 6 through 21, as well as the number of $f$-evaluations, is $O \left( \frac{k^2 \cdot \log^2 N \cdot \log \left(\frac{k \cdot \ln N}{\epsilon} \right) \cdot \log^2 \left(\frac{\epsilon N}{k}\right)}{\log^2 \left(\frac{k}{\epsilon} \right) \cdot \epsilon^2 \cdot \log \log \left(\frac{\epsilon N}{k}\right)} \right)$.
\label{thm:Recon2}
\end{thm}

\noindent \textit{Proof:}  See Appendix~\ref{app:proof_thm:Recon2}.~~$\Box$\\

The overall runtime behavior of Algorithm~\ref{alg:reconstruct2} is determined by both the runtime of Algorithm~\ref{alg:fastMult} as called in line 4 of Algorithm~\ref{alg:reconstruct2}, and the runtime stated in Theorem~\ref{thm:Recon2}.  The overall runtime complexity of Algorithm~\ref{alg:reconstruct2} is therefore given in Equation~\ref{eqn:Alg1RUNTIME2}.  As in Section~\ref{sec:FourierRecov} above, both this runtime and the number of function evaluations required for approximate Fourier reconstruction can be decreased by reducing the number of measurement matrix rows (i.e., $\mathcal{R}_{\lambda,K}$ rows) used to estimate each Fourier coefficient.  This effectively replaces $K$ in Algorithm~\ref{alg:reconstruct2} with a significantly smaller value (e.g., the value $l$ from Corollary~\ref{cor:RandSj}).  However, in exchange for the resulting runtime improvements we will once again have to sacrifice approximation guarantees for a small probability of outputting a highly inaccurate answer.

Following the strategy above, we will improve the performance of Algorithm~\ref{alg:reconstruct2} by modifying its utilized measurement matrix as follows:  Instead of using a $\mathcal{M}_{s_1,K}$ matrix as constructed in Section~\ref{sec:makeM} to build $\mathcal{R}_{\lambda,K} = \mathcal{M}_{s_1,K} \circledast \mathcal{N}_{\lambda,s_1}$, we will instead use a randomly constructed $\mathcal{M}_{\tilde{S}}$ matrix as described in Section~\ref{sec:FixedSig} to build $\mathcal{R}_{\lambda,\tilde{S}} = \mathcal{M}_{\tilde{S}} \circledast \mathcal{N}_{\lambda,s_1}$.  Corollary~\ref{cor:RandSj} combined with the proof of Lemma~\ref{lem:Rprops} ensures that such a randomly constructed measurement matrix, $\mathcal{R}_{\lambda,\tilde{S}}$, will be likely to have all the properties of $\mathcal{R}_{\lambda,K}$ matrices that Algorithm~\ref{alg:reconstruct2} needs to function correctly.  Hence, with high probability we will receive output from Algorithm~\ref{alg:reconstruct2} with the same approximation error bounds as derived for Theorem~\ref{thm:Recon2}.  Formalizing these ideas we obtain the following Corollary proven in Appendix~\ref{app:proof_cor:Recon2}.

\begin{cor}
Suppose $f: [0,2\pi] \rightarrow \mathbbm{C}$ has $\hat{f} \in l^1$. Let $\sigma \in [2/3,1)$ and $N, k, \epsilon^{-1} \in \mathbbm{N}-\{ 1 \}$ with $N > (k/\epsilon) \geq 2$.  Algorithm~\ref{alg:reconstruct2} may be executed using a random matrix, $\mathcal{R}_{\lambda,\tilde{S}} = \mathcal{M}_{\tilde{S}} \circledast \mathcal{N}_{\lambda,s_1}$, in place of the deterministic matrix, $\mathcal{R}_{\lambda,K} = \mathcal{M}_{s_1,K} \circledast \mathcal{N}_{\lambda,s_1}$, considered above.  In this case Algorithm~\ref{alg:reconstruct2} will produce an output vector, $\vec{x}_{S} \in \mathbbm{C}^N$, that satisfies Equation~\ref{eqn:Recon2Error} with probability at least $\sigma$.  Both the runtime of lines 6 through 21 and the number of points in $[0,2\pi]$ at which $f$ will be evaluated are 
$$O \left( \frac{k}{\epsilon} \cdot \log^3 N \cdot \log \left( \frac{N}{1 - \sigma} \right) \right).$$
Explicit upper bounds on the number of point evaluations are easily obtained from the proof below.
\label{cor:Recon2}
\end{cor}

\noindent \textit{Proof:}  See Appendix~\ref{app:proof_cor:Recon2}.~~$\Box$\\

When executed with a random matrix, $\mathcal{R}_{\lambda,\tilde{S}}$, as input the overall runtime complexity of Algorithm~\ref{alg:reconstruct2} will be determined by both the runtime stated in Corollary~\ref{cor:Recon2} and the runtime of Algorithm~\ref{alg:fastMult}.  Suppose $\tilde{S}$ is a subset of $O\left( \log \left( \frac{N}{1 - \sigma} \right) \right)$ $s_j$ values defined as per Equations~\ref{eqn:Primes} --~\ref{eqn:Def_k}.  Then, Algorithm~\ref{alg:fastMult} will have a runtime complexity of
\begin{align}
O\left( \sum^{\lambda}_{i=1} \sum_{s_j \in \tilde{S}} t_i \cdot s_{j} \log s_{j} \right) & = O\left( \frac{p_{q+K} \cdot \log p_{q+K} \cdot p^2_{\lambda}}{\log p_{\lambda}} \cdot \log \left( \frac{N}{1 - \sigma} \right) \right) ~~\textrm{  (see~~\cite{IDSFA}, Corollary~\ref{cor:RandSj})} \nonumber \\ & = 
O\left( \frac{k \cdot \log_{(k/\epsilon)} N \cdot \log^2 \left( \frac{k \cdot \log N }{\epsilon} \right) \cdot \ln^2 \left( \frac{\epsilon \cdot N}{k} \right)}{\epsilon \cdot \ln \ln \left( \frac{\epsilon \cdot N}{k}\right)} \cdot \log \left( \frac{N}{1 - \sigma} \right) \right) ~~\textrm{  (see~Equation~\ref{equ:pq+K},~Lemma~\ref{lem:RowBound})}.
\label{eqn:Alg2RandRUNTIME}
\end{align}
Thus, if we are willing to fail with probability at most $1 - \sigma = 1 / N^{O(1)}$, then Algorithm~\ref{alg:reconstruct2} executed with a random input matrix will have a total runtime complexity of $O\left( (k/\epsilon) \cdot \log^4 N \cdot \log \left( \frac{k \cdot \log N }{\epsilon} \right) \right)$.  

\section{Higher Dimensional Fourier Transforms}
\label{sec:HighDim}

In this section we will consider methods for approximating the Fourier transform of a periodic function of $D$ variables, $f: [0,2\pi]^D \rightarrow \mathbbm{C}$.  To begin, we will demonstrate how to approximate the Fourier transform of $f$ by calculating the discrete Fourier transform of a related one-dimensional function, $f_{\rm new}: [0, 2\pi] \rightarrow \mathbbm{C}$.  This dimensionality reduction technique for multidimensional Fourier transforms will ultimately enable us to quickly approximate $\hat{f}$ by applying the methods of Section~\ref{sec:FourierRecov2} to $f$'s related one dimensional function $f_{\rm new}$.  The end result will be a set of algorithms for approximating $\hat{f}$ whose runtimes scale polynomially in the input dimension $D$.

Suppose that the Fourier transform of $f$ above, $\hat{f}: \mathbbm{Z}^D \rightarrow \mathbbm{C}$, is near zero for all integer points outside of the $D$-dimensional cubic lattice $\left( [-M/2, M/2] \cap \mathbbm{Z} \right)^D$.  In order to help us approximately recover $\hat{f}$ we will choose $D$ pairwise relatively prime integers, $P_1, \dots, P_D \in \mathbbm{N}$, with the property that $P_d > M \cdot D$ for all $d \in [1,D] \cap \mathbbm{N}$.  Set $\tilde{N} = \prod^D_{d=1} P_d$.  Furthermore, let $y^{-1 \textrm{ mod } p} \in [0, p) \cap \mathbbm{N}$ denote the multiplicative inverse of $\left( y \textrm{ mod } p \right) \in \mathbbm{Z}_p$ when it exists.  Note that $y^{-1 \textrm{ mod } p}$ will exist whenever $y$ is relatively prime to $p$.  

We may now define the function $f_{\rm new}: [0, 2\pi] \rightarrow \mathbbm{C}$ to be
\begin{equation}
f_{\rm new}(x) = f \left( \frac{\tilde{N}}{P_1}x, \frac{\tilde{N}}{P_2}x, \dots, \frac{\tilde{N}}{P_D}x \right).
\label{def:fnew}
\end{equation}
Considering the Fourier transform of $f_{\rm new}$ we can see that
\begin{align}
\hat{f}_{\rm new}(\omega) &~=~ \frac{1}{2 \pi} \int^{2\pi}_{0} \mathbbm{e}^{- \mathbbm{i} \omega x} ~f_{\rm new} (x)~ dx ~=~ \frac{1}{2 \pi} \sum_{\left( \omega_1, \dots, \omega_D \right) \in \mathbbm{Z}^D} \hat{f}\left( \omega_1, \dots, \omega_D \right) ~\int^{2\pi}_{0} \mathbbm{e}^{-\mathbbm{i} x \left( \omega - \sum^D_{d=1} \frac{\tilde{N}}{P_d}\omega_d \right)} ~dx\nonumber \\ &~=~ \sum_{\left( \omega_1, \dots, \omega_D \right) \in \mathbbm{Z}^D~\textrm{s.t.}~\omega=\sum^D_{d=1} \frac{\tilde{N}}{P_d}\omega_d} \hat{f}\left( \omega_1, \dots, \omega_D \right). \label{eqn:fnewFourier}
\end{align}
Recall that we are primarily interested in capturing the information about $\hat{f}$ inside $\left( [-M/2, M/2] \cap \mathbbm{Z} \right)^D$.  Looking at the $\omega \in \mathbbm{Z}$ for which $\hat{f}_{\rm new}$ can impacted by $\left( \omega_1, \dots, \omega_D \right) \in \left( [-M/2, M/2] \cap \mathbbm{Z} \right)^D$ we can see that
$$\left| \omega \right| ~\leq~ \sum^D_{d=1} \left| \frac{\omega_d \tilde{N}}{P_d} \right| ~\leq~ \sum^D_{d=1}  \frac{M \tilde{N}}{2 P_d}  ~<~ \sum^D_{d=1} \frac{\tilde{N}}{2 D}  ~=~ \frac{\tilde{N}}{2}.$$
Hence, we may consider $f_{\rm new}$ to have an effective bandwidth of $\tilde{N}$.

More importantly, there is a bijective correspondence between the integer lattice points, $\left( \omega_1, \dots, \omega_D \right) \in \left( [-M/2, M/2] \cap \mathbbm{Z} \right)^D$, and their representative frequency, $ \omega \in [-\tilde{N}/2, \tilde{N}/2] \cap \mathbbm{Z}$, in $\hat{f}_{\rm new}$.  Define the function 
$$g: \left(-\frac{P_1}{2}, \frac{P_1}{2} \right] \cap \mathbbm{N} \times \dots \times \left(-\frac{P_D}{2}, \frac{P_D}{2} \right] \cap \mathbbm{N} \rightarrow \left(-\frac{\tilde{N}}{2}, \frac{\tilde{N}}{2} \right] \cap \mathbbm{N}$$
to be 
$$g(x_{1}, \dots, x_{D}) = \left( \sum^{D}_{d = 1} \left(\frac{\tilde{N}}{P_d} \right) \cdot x_{d} \right) \textrm{ mod } \tilde{N}.$$
The Chinese Remainder Theorem tells us that $g$ is a well-defined bijection.  Furthermore, it is not difficult to see that 
$$g^{-1}(x) = \left( x \cdot \left(\tilde{N}/P_1 \right)^{-1 \textrm{ mod } P_1} \textrm{ mod } P_1,~\dots,~ x \cdot \left(\tilde{N}/P_D \right)^{-1 \textrm{ mod } P_D} \textrm{ mod } P_D \right).$$
Thus, we have $\hat{f}_{\rm new}(\omega) \approx \hat{f}\left( g^{-1}(\omega) \right)$.

We now have a three-step algorithm for finding a sparse Fourier approximation for any function $f: [0,2\pi]^D \rightarrow \mathbbm{C}$.  All we must do is: $(i)$ Implicitly create $f_{\rm new}$ as per Equation~\ref{def:fnew}, $(ii)$ Use the techniques from Section~\ref{sec:FourierRecov2} to approximate $\vec{\hat{f}}_{\rm new}$, and then $(iii)$ Use the approximation for $\vec{\hat{f}}_{\rm new}$ to approximate $\vec{\hat{f}}$ via Equation~\ref{eqn:fnewFourier}.  The following theorem summarizes some of the results one can achieve by utilizing this approach.

\begin{thm}
Suppose $f: [0,2\pi]^D \rightarrow \mathbbm{C}$ is bandlimited so that $\hat{f}(\omega_1, \dots, \omega_D) = 0$ if $\left( \omega_1, \dots, \omega_D \right) \notin$ $\left( [-\frac{M}{2}, \frac{M}{2}] \cap \mathbbm{Z} \right)^D$.  Define $\tilde{N}$ as above and suppose that $\tilde{N}, k, \epsilon^{-1} \in \mathbbm{N}-\{ 1 \}$ with $\tilde{N} > (k/\epsilon)^2 \geq 4$.  Then, Algorithm~\ref{alg:reconstruct2} combined with the bijective mapping, $g$, above will output an $\vec{x}_{S} \in \mathbbm{C}^{\tilde{N}}$ satisfying
\begin{equation}
\left\| \hat{f} - \left(\vec{x}_{S} \circ g\right) \right\|_2 \leq \left\| \hat{f} - \hat{f}^{\rm~ opt}_k \right\|_2 + \frac{22\epsilon \cdot \left\| \hat{f} - \hat{f}^{\rm~ opt}_{(k/\epsilon)} \right\|_1}{\sqrt{k}}.
\label{eqn:MultReconError}
\end{equation} 
Both the runtime of lines 6 -- 21, and the number of points in $[0,2\pi]^D$ at which $f$ will be evaluated, will be 
$$O \left( \frac{k^2 \cdot D^4 \cdot \log^4 (MD)}{\log \left(\frac{k}{\epsilon} \right) \cdot \epsilon^2} \right).$$

If succeeding with probability $\sigma \in [2/3,1)$ is sufficient, and $\tilde{N} > (k/\epsilon) \geq 2$, Algorithm~\ref{alg:reconstruct2} may instead be executed using a random matrix $\mathcal{R}_{\lambda,\tilde{S}}$.  In this case Algorithm~\ref{alg:reconstruct2} will produce an output vector, $\vec{x}_{S} \in \mathbbm{C}^{\tilde{N}}$, that satisfies Equation~\ref{eqn:MultReconError} with probability at least $\sigma$.  Both the runtime of lines 6 -- 21, and the number of points in $[0,2\pi]^D$ at which $f$ will be evaluated, will be 
$$O \left( \frac{k \cdot D^4}{\epsilon} \cdot \log^3 (MD) \cdot \log \left( \frac{MD}{1 - \sigma} \right) \right).$$

Finally, if an exponential runtime of $\Omega \left( (DM)^D \right)$ is acceptable, we note that both Theorem~\ref{thm:Recon1} and Corollary~\ref{cor:RandRecon1} can also be adapted to recovering $f: [0,2\pi]^D \rightarrow \mathbbm{C}$ by substituting $N$ with $\Theta\left( (MD)^D \right)$ everywhere in their statements.
\label{thm:MultDimRecov}
\end{thm}

\noindent \textit{Proof:}  See Appendix~\ref{app:MultDimRecov}.~~$\Box$\\

Note that traditional FFT algorithms (e.g., \cite{FFT,rabiner-schafer-rader,1968-bluestein}) require $\Omega\left(M^D \right)$-time to calculate the Fourier transform of a bandlimited function $f: [0,2\pi]^D \rightarrow \mathbbm{C}$.  In contrast, Theorem~\ref{thm:MultDimRecov} allows us to approximate $\hat{f}$ using exponentially fewer (in $D$) operations.  Hence, if $f$ has a relatively sparse Fourier representation (e.g., if $\hat{f}$ is dominated by $k = M^{o(D)}$ energetic frequencies), Theorem~\ref{thm:MultDimRecov} allows $\hat{f}$ to be accurately approximated much more quickly than possible using standard techniques.

\section{Conclusion}
\label{sec:Conc}

In conclusion, it is worth pointing out that the methods developed in this paper for approximating the Fourier transforms of periodic functions are also applicable to the approximation of functions which have accurate sparse representations in related bases. For example, all the theorems proven herein will also apply to functions with sparsely representable Cosine or Chebyshev expansions (see \cite{BoydAl} for an in depth discussion of the relationships between these series expansions).  Hence, we have also implicitly constructed sublinear-time algorithms for approximating these related transforms.

\section*{Acknowledgments}

I would like to thank Martin Strauss and Holger Rauhut for answering questions about their work.  I would also like to thank Martin Strauss, Holger Rauhut, and Tsvetanka Sendova for helpful discussions and advice.  The vast majority of this paper was written while the author was supported as a postdoctoral fellow at the Institute for Mathematics and its Applications (IMA).

\bibliographystyle{abbrv}
\bibliography{JournalFFT}

\begin{thebibliography}{10}

\bibitem{1968-bluestein}
L.~I. Bluestein.
\newblock {A Linear Filtering Approach to the Computation of Discrete Fourier
  Transform}.
\newblock {\em IEEE Transactions on Audio and Electroacoustics}, 18:451--455,
  1970.

\bibitem{HardThreshforCS}
T.~Blumensath and M.~E. Davies.
\newblock Iterative hard thresholding for compressed sensing.
\newblock {\em Applied and Computational Harmonic Analysis}, 27(3):265 -- 274,
  2009.

\bibitem{BoydAl}
J.~P. {Boyd}.
\newblock {\em {Chebyshev and Fourier Spectral Methods}}.
\newblock Dover Publications, Inc., 2001.

\bibitem{CS1}
E.~Candes, J.~Romberg, and T.~Tao.
\newblock Robust uncertainty principles: Exact signal reconstruction from
  highly incomplete frequency information.
\newblock {\em IEEE Trans. Inform. Theory}, 52:489--509, 2006.

\bibitem{CS4}
E.~Candes, J.~Romberg, and T.~Tao.
\newblock Stable signal recovery from incomplete and inaccurate measurements.
\newblock {\em Communications on Pure and Applied Mathematics},
  59(8):1207--1223, 2006.

\bibitem{NearOpt}
E.~Candes and T.~Tao.
\newblock Near optimal signal recovery from random projections: Universal
  encoding strategies?
\newblock {\em IEEE Trans. on Information Theory}, 2006.

\bibitem{BestkTerm}
A.~Cohen, W.~Dahmen, and R.~DeVore.
\newblock {Compressed Sensing and Best $k$-term Approximation}.
\newblock {\em Journal of the American Mathematical Society}, 22(1):211--231,
  January 2008.

\bibitem{FFT}
J.~Cooley and J.~Tukey.
\newblock An algorithm for the machine calculation of complex {F}ourier series.
\newblock {\em Math. Comput.}, 19:297--301, 1965.

\bibitem{MedianMedian}
D.~Dor and U.~Zwick.
\newblock Median selection requires (2+eps)n comparisons.
\newblock {\em SIAM Journal of Discrete Mathematics}, 14:125--134, 2001.

\bibitem{PrimeBound}
P.~Dusart.
\newblock The $k^{\rm th}$ prime is greater than $k(\ln k + \ln\ln k-1)$ for $k
  \geq 2$.
\newblock {\em Mathematics of Computation}, 68(225), 1999.

\bibitem{HolgerCSBook}
M.~Fornasier and H.~Rauhut.
\newblock Compressive sensing.
\newblock {\em Handbook of Mathematical Methods in Imaging}, Springer, to
  appear.

\bibitem{AAFFT1}
A.~Gilbert, S.~Guha, P.~Indyk, S.~Muthukrishnan, and M.~Strauss.
\newblock Near-optimal sparse {F}ourier estimation via sampling.
\newblock {\em ACM STOC}, pages 152--161, 2002.

\bibitem{AAFFT2}
A.~Gilbert, S.~Muthukrishnan, and M.~Strauss.
\newblock Improved time bounds for near-optimal sparse {F}ourier
  representations.
\newblock {\em Proceedings of SPIE Wavelets XI}, 2005.

\bibitem{JournalDSFT}
M.~A. Iwen.
\newblock Combinatorial sublinear-time fourier algorithms.
\newblock {\em Foundations of Computational Mathematics}, 10(3):303 -- 338,
  2010.

\bibitem{IDSFA}
M.~A. Iwen and C.~V. Spencer.
\newblock Improved bounds for a deterministic sublinear-time sparse fourier
  algorithm.
\newblock In {\em Conference on Information Sciences and Systems (CISS)}, 2008.

\bibitem{SortSearch}
D.~Knuth.
\newblock {\em {The Art of Computer Programming. Volume 3: Sorting and
  Searching}}.
\newblock Addison-Wesley, 1998.

\bibitem{OMPvsBP}
S.~Kunis and H.~Rauhut.
\newblock {Random Sampling of Sparse Trigonometric Polynomials II - Orthogonal
  Matching Pursuit versus Basis Pursuit}.
\newblock {\em Foundations of Computational Mathematics}, 8(6):737--763, 2008.

\bibitem{RandAlg}
R.~Motwani and P.~Raghavan.
\newblock {\em Randomized Algorithms}.
\newblock Cambridge University Press, 1995.

\bibitem{COSAMP}
D.~Needell and J.~A. Tropp.
\newblock Cosamp: Iterative signal recovery from incomplete and inaccurate
  samples.
\newblock {\em Applied and Computational Harmonic Analysis}, 26(3), 2008.

\bibitem{ROMP}
D.~Needell and R.~Vershynin.
\newblock Uniform uncertainty principle and signal recovery via regularized
  orthogonal matching pursuit.
\newblock {\em Foundations of Computational Mathematics}, 9:317--334, 2009.

\bibitem{ROMPstable}
D.~Needell and R.~Vershynin.
\newblock Signal recovery from incomplete and inaccurate measurements via
  regularized orthogonal matching pursuit.
\newblock {\em IEEE Journal of Selected Topics in Signal Processing}, to
  appear.

\bibitem{NumTheory}
I.~Niven, H.~S. Zuckerman, and H.~L. Montgomery.
\newblock {\em {An Introduction to The Theory of Numbers}}.
\newblock John Wiley \& Sons, Inc., 1991.

\bibitem{rabiner-schafer-rader}
L.~Rabiner, R.~Schafer, and C.~Rader.
\newblock {T}he {C}hirp z-{T}ransform {A}lgorithm.
\newblock {\em IEEE Transactions on Audio and Electroacoustics},
  AU-17(2):86--92, June 1969.

\bibitem{HolgerStableRecov}
H.~Rauhut.
\newblock Compressive sensing and structured random matrices.
\newblock {\em Theoretical Foundations and Numerical Methods for Sparse
  Recovery}, vol. 9 in Radon Series Comp. Appl. Math.:1--92, deGruyter, 2010.

\bibitem{NumThyBound2}
G.~Robin.
\newblock Estimation de la fonction de tchebychef {$\theta$} sur le k-ieme
  nombre premier et grandes valeurs de la fonctions {$\omega(n)$}, nombre de
  diviseurs premiers de n.
\newblock {\em Acta Arithmetica}, 42(4):367 -- 389, 1983.

\bibitem{CSb1}
M.~Rudelson and R.~Vershynin.
\newblock Sparse reconstruction by convex relaxation: Fourier and gaussian
  measurements.
\newblock In {\em 40th Annual Conference on Information Sciences and Systems
  (CISS)}, 2006.

\bibitem{CS2}
J.~Tropp and A.~Gilbert.
\newblock Signal recovery from partial information via orthogonal matching
  pursuit.
\newblock {\em IEEE Trans. Info. Theory}, 53(12):4655--4666, Dec. 2007.

\end{thebibliography}

\appendix

\section{Proof of Theorem~\ref{thm:RowBound}}
\label{app:proof_thm:RowBound}

Let $\pi(n)$ be the number of primes no greater than $n$.  In \cite{PrimeBound} it is shown that
$$\frac{n}{\ln n} \left(1 + \frac{0.992}{\ln n} \right) \leq \pi(n) \leq \frac{n}{\ln n} \left(1 + \frac{1.2762}{\ln n} \right)$$
for all $n \geq 599$.  Using this result (in combination with numerical tests for $n < 600$) we obtain the following bounds for $q+K$ and $q$ (see Equation~\ref{eqn:Def_k}).
\begin{equation}
q + K \leq \pi(k/\epsilon) + K + 1 \leq \frac{k \left\lfloor \log_{(k/\epsilon)} N \right\rfloor }{\epsilon} \left( c + \frac{1}{\ln (k/\epsilon) \cdot \left\lfloor \log_{(k/\epsilon)} N \right\rfloor} + \frac{1.2762}{\ln^2 (k/\epsilon) \cdot \left\lfloor \log_{(k/\epsilon)} N \right\rfloor} + \frac{2 \cdot \epsilon}{k \cdot \left\lfloor \log_{(k/\epsilon)} N \right\rfloor} \right).
\label{eqn:q+Kbound}
\end{equation}
and 
\begin{equation}
q \geq \pi(k/\epsilon) \geq \max \left\{ \frac{k}{\epsilon \cdot \ln (k/\epsilon)} \left(1 + \frac{0.992}{\ln (k/\epsilon)} - \frac{8.85 \cdot \epsilon}{k} \right), 1 \right\}.
\label{equ:qBound}
\end{equation}
Continuing, we can bound $m$ if our $s_j$ values are chosen to be primes as above by noting that
\begin{align}
\sum^{q-1}_{j=1} p_j & \geq \sum^{q-1}_{j=1} j \cdot \ln(j) ~~\textrm{  (see \cite{PrimeBound}) } \nonumber \\ & \geq \int^{q-1}_{1} x \cdot \ln x ~dx \geq \frac{(q-1)^2}{2} \left( \ln(q-1) - \frac{1}{2} \right)
\label{eqn:LowerSum}
\end{align}
and
\begin{align}
\sum^{q+K-1}_{j=1} p_j & \leq 10 + \sum^{q+K-1}_{j=4} j \cdot \ln(p_j) ~~\textrm{  (see \cite{PrimeBound}) } \nonumber \\ & \leq 10 + \ln(p_{q+K}) \cdot \left(\sum^{q+K-1}_{j=4} j \right) \leq \frac{(q+K-1)(q+K)}{2} \cdot \ln \big( (q+K)\cdot \left( \ln(q+K) + \ln \ln (q+K) \right) \big) ~~\textrm{  (see \cite{PrimeBound}) } \label{eqn:UpperSum}
\\ & \leq \frac{3}{4}(q+K)^2 \cdot \ln(q+K). 
\label{eqn:upperSum}
\end{align}
Using Equation~\ref{eqn:q+Kbound} together with Equation~\ref{eqn:upperSum} finishes the proof.  More specifically, we have that 
$$\left( c + \frac{1}{\ln (k/\epsilon) \cdot \left\lfloor \log_{(k/\epsilon)} N \right\rfloor} + \frac{1.2762}{\ln^2 (k/\epsilon) \cdot \left\lfloor \log_{(k/\epsilon)} N \right\rfloor} + \frac{2 \cdot \epsilon}{k \cdot \left\lfloor \log_{(k/\epsilon)} N \right\rfloor} \right) \leq \left(c + \frac{1}{\ln 4} + \frac{1.2762}{\ln^2 4} + \frac{1}{2} \right) \leq (c + 1.89).$$
Therefore, we can see that
$$m \leq \sum^{q+K-1}_{j=1} p_j \leq \frac{3 (c + 1.89)^2 \cdot k^2 \left\lfloor \log_{(k/\epsilon)} N \right\rfloor^2}{4 \cdot \epsilon^2} \cdot \ln\left( \frac{(c + 1.89) \cdot k \left\lfloor \log_{(k/\epsilon)} N \right\rfloor}{\epsilon} \right)$$
as we wished to prove.

\section{Proof of Corollary~\ref{cor:RandSj}}
\label{app:proof_cor:RandSj}

We prove the result via an argument similar to the one used to prove Lemma~2 in \cite{JournalDSFT}.  Fix $n \in S$.  We will select our multiset of $s_j$ values, $\tilde{S}$, by independently choosing $l$ elements of $\{ s_1, s_2, \dots, s_K \}$ uniformly at random with replacement.  The first element chosen for $\tilde{S}$ will be denoted $s_{j_1}$, the second $s_{j_2}$, and so forth.  Let $Q^n_h$ be the random variable indicating whether the $s_{j_h}$ value selected for $\tilde{S}$ satisfies
\begin{equation}
\left| \left(\mathcal{M}_{s_1,K,n} \vec{x}\right)_{j_h} - x_n \right| \leq \frac{ \epsilon \cdot \left\| \vec{x} - \vec{x}^{\rm opt}_{(k/\epsilon)} \right\|_1}{k}.
\label{eqn:PropRand}
\end{equation}
Therefore, 
$$Q^n_h = \left\{ \begin{array}{ll} 1 & \textrm{if } s_{j_h} \textrm{ satisfies Property~\ref{eqn:PropRand}} \\ 0 & {\rm otherwise} \end{array} \right..$$
Theorem~\ref{thm:Mest} tells us that $\mathbbm{P}\left[ Q^n_h = 1 \right] > \frac{6}{7}$.  Furthermore, $\mu = \mathbbm{E} \left[ \sum^{l}_{h = 1} Q^n_h \right] \geq \frac{6 \cdot l}{7}$.

Using the Chernoff bound (see \cite{RandAlg}) we get that the probability of $$\sum^{l}_{h = 1} Q^n_h < \frac{4 \cdot l}{7}$$
is less than $e^{- \frac{\mu}{18}} \leq e^{-\frac{l}{21}} \leq \frac{1 - \sigma}{|S|}$.  Since $l > 21$ we can see that $\sum^{l}_{h = 1} Q^n_h$ will be less than $\frac{l + 1}{2}$ with probability less than $\frac{1-\sigma}{|S|}$.  Hence, Property~\ref{eqn:PropRand} will be satisfied by more than $l/2$ of the $s_{j_h} \in \tilde{S}$ with high probability.  Applying the union bound shows that the majority of the entries in $\tilde{S}$ will indeed satisfy Property~\ref{eqn:PropRand} for all $n \in S$ with probability at least $\sigma$.  The result follows.

\section{Proof of Corollary~\ref{cor:RandRecon1}}
\label{app:proof_cor:RandRecon1}

Apply Corollary~\ref{cor:RandSj} with $c = 14$, $\vec{x} = \vec{\hat{f}}$, and $S = \left( - \left\lceil \frac{N}{2} \right\rceil, \left\lfloor \frac{N}{2} \right\rfloor \right] \cap \mathbbm{Z}$ to obtain $\tilde{S}$, a multiset of $\left\lceil 21 \cdot \ln \left( \frac{N}{1 - \sigma} \right) \right\rceil$ $s_j$ values.  With probability at least $\sigma$ more than half (with multiplicity) of the entries of $\mathcal{M}_{\tilde{S}},\omega \vec{\hat{f}}$ will estimate $\hat{f}_{\omega}$ to within $(\epsilon/k) \cdot \left\| \vec{\hat{f}} - \vec{\hat{f}}^{\rm ~opt}_{(k/\epsilon)} \right\|_1$ precision for all $\omega \in \left( - \left\lceil \frac{N}{2} \right\rceil, \left\lfloor \frac{N}{2} \right\rfloor \right] \cap \mathbbm{Z}$.  Furthermore, $\mathcal{M}_{\tilde{S}} \vec{\hat{f}}$ can still be approximately computed using Algorithm~\ref{alg:fastMult} if only the unique $s_j$ values in $\tilde{S}$ are given as the relatively prime inputs.  In this case Lemma~\ref{lem:MultError} will also still hold.  Taken all together we can see that with probability at least $\sigma$ all $N$ $x_{\omega}$ values produced by lines 6 and 7 of Algorithm~\ref{alg:reconstruct1} will have 
$$\left| x_{\omega} - \hat{f}_\omega \right| \leq \sqrt{2} \cdot \left(\frac{ \epsilon \cdot \left\| \vec{\hat{f}} - \vec{\hat{f}}^{\rm ~opt}_{(k/\epsilon)} \right\|_1}{k} + \left\| ~\hat{f} - \bar{\hat{f}} ~ \right\|_1 \right).$$
The Equation~\ref{eqn:Recon1Error} error bound now follows from the proof of Theorem~\ref{thm:Recon1}.

To upper bound the number of required function evaluations we will bound the number of rows for a particular $\mathcal{M}_{\tilde{S}}$ matrix constructed with primes as per Section~\ref{sec:makeM}.  In particular, we will assume that 
$\tilde{S}$ contains at most $\left\lceil 21 \cdot \ln \left( \frac{N}{1 - \sigma} \right) \right\rceil$ individual $s_j$ values defined as in Equations~\ref{eqn:Primes} --~\ref{eqn:Def_k} with $K = 14 \cdot (k / \epsilon) \left \lfloor \log_{s_1} N \right \rfloor + 1$.  In this case Equation~\ref{eqn:q+Kbound} together with results from \cite{PrimeBound} tell us that $s_K$ is at most 
\begin{equation}
15.89 \cdot \frac{k \left\lfloor \log_{(k/\epsilon)} N \right\rfloor }{\epsilon} \cdot \left( \ln \left( \frac{15.89 \cdot k \left\lfloor \log_{(k/\epsilon)} N \right\rfloor }{\epsilon} \right) + \ln \ln \left( \frac{15.89 \cdot k \left\lfloor \log_{(k/\epsilon)} N \right\rfloor }{\epsilon} \right) \right).
\label{eqn:SKboundApp}
\end{equation}
The stated upper bound on the number of required function evaluations follows.  The stated runtime follows from the fact that each line 6 and 7 median now only involves $O\left( \log \left( \frac{N}{1 - \sigma} \right) \right)$ values.

\section{Proof of Lemma~\ref{lem:RowBound}}
\label{app:proof_lem:RowBound}

We can always set $t_1 = p_1 < \cdots < t_\lambda  = p_\lambda$.  In this case we require that $p_\lambda < s_1 \leq$ the smallest prime factor of $s_1, \dots, s_K$.  Secondly, we require that $\sum^\lambda_{i=1} \ln p_i \geq \ln \left( \frac{N}{s_1} \right).$  Using results from \cite{NumThyBound2} it is easily verified that 
$$\sum^\lambda_{i=1} \ln p_i \geq \lambda \cdot \left( \ln \lambda - 1 \right)$$  
for all $\lambda \in \mathbbm{N}^{+}$.  Setting $\lambda = \left\lceil 3 \ln \left( \frac{N}{s_1} \right) / \ln \ln \left( \frac{N}{s_1} \right) \right\rceil$ in the equation above we can see that
$$\sum^\lambda_{i=1} \ln p_i \geq \ln \left( \frac{N}{s_1} \right) \cdot 3 \left( 1 - \frac{\ln \ln \ln \left( \frac{N}{s_1} \right) }{\ln \ln \left( \frac{N}{s_1} \right) } \right) \geq \ln \left( \frac{N}{s_1} \right)$$ 
as long as $N/s_1 \geq 3$.  Hence, if we choose our $t_i$ values to be the first $\lambda$ primes the second requirement will be satisfied.

Results from \cite{PrimeBound} then tell us that
$$t_\lambda = p_\lambda \leq p_{\left\lceil 3 \ln \left( N / s_1 \right) / \ln \ln \left( N / s_1 \right) \right\rceil} \leq \left\lceil 3 \cdot \frac{\ln \left( N / s_1 \right)}{\ln \ln \left( N / s_1 \right)} \right\rceil \cdot \left( \ln \left\lceil 3 \cdot \frac{\ln \left( N / s_1 \right)}{\ln \ln \left( N / s_1 \right)} \right\rceil + \ln \ln \left\lceil 3 \cdot \frac{\ln \left( N / s_1 \right)}{\ln \ln \left( N / s_1 \right)} \right\rceil \right) < s_1.$$
Therefore, the prime $t_i$ values we have selected will also satisfy the first requirement above.  To bound the smallest possible number of rows we note that
$$\tilde{m} \leq 1 + \sum^{\left\lceil 3 \ln \left( \frac{N}{s_1} \right) / \ln \ln \left( \frac{N}{s_1} \right) \right\rceil}_{i = 1} p_i ~\leq~ \frac{3}{4}\left(\left\lceil 3 \cdot \frac{\ln \left( N / s_1 \right)}{\ln \ln \left( N / s_1 \right)} \right\rceil + 1\right)^2 \cdot \ln\left(\left\lceil 3 \cdot \frac{\ln \left( N / s_1 \right)}{\ln \ln \left( N / s_1 \right)} \right\rceil+1 \right) +1~~\textrm{(see Equation~\ref{eqn:upperSum})}.$$
The stated result follows.

\section{Proof of Lemma~\ref{lem:Rprops}}
\label{app:proof_lem:Rprops}

In addition to $\vec{x}$ we will also consider $\vec{y} \in \mathbbm{C}^N$ defined by
$$y_{n'} = \left| x_{n'} \right|~\textrm{for all}~n' \in [0,N) \cap \mathbbm{N}.$$
Note that $\vec{y}$ and $\vec{x}$ will not only share the same optimal $(k/ \epsilon)$-term support subset, $S^{\rm opt}_{(k/ \epsilon)} \subset [0,N) \cap \mathbbm{N}$, but will also have $\left\| \vec{x} - \vec{x}^{\rm opt}_{(k/\epsilon)} \right\|_1 = \left\| \vec{y} - \vec{y}^{\rm opt}_{(k/\epsilon)} \right\|_1$.  Theorem~\ref{thm:Mest} tells us that more than $\frac{K}{2}$ entries of $\mathcal{M}_{s_1,K,n} \cdot \vec{y}$ will estimate $y_n$ to within $\frac{\epsilon \cdot \left\| \vec{y} - \vec{y}^{\rm opt}_{(k/\epsilon)} \right\|_1}{k} = \bar{\delta} = \frac{\epsilon \cdot \left\| \vec{x} - \vec{x}^{\rm opt}_{(k/\epsilon)} \right\|_1}{k}$ precision.  Let $\left( \mathcal{M}_{s_1,K,n} \cdot \vec{y} \right)_{j'}$ for $j' \in [1,K] \cap \mathbbm{N}$ be one of these $\frac{K}{2}$ entries.  The proof of Lemma~\ref{lem2} tells us that the row associated with this entry also has the property that
$$\left| \left( \mathcal{M}_{s_1,K,n} \cdot \vec{x} \right)_{j'} - x_{n} \right| \leq \sum_{n' \equiv n \textrm{~mod~} s_{j'},~n' \notin S^{\rm opt}_{(k / \epsilon)}, ~n' \neq n} y_{n'} = \left| \left( \mathcal{M}_{s_1,K,n} \cdot \vec{y} \right)_{j'} - y_{n} \right| \leq \bar{\delta}.$$ 
Therefore, we have established property $(1)$.

Considering property $(2)$ for this $j'$ we can see that for all $i \in [1,\lambda] \cap \mathbbm{N}$ we will have
\begin{align}
\left| \bar{r}_{i,j',n \textrm{~mod~} t_i \cdot s_{j'}} \cdot \vec{x} - x_{n} \right| = \left| \sum_{n' \equiv n \textrm{~mod~} t_i \cdot s_{j'},~n' \notin S^{\rm opt}_{(k / \epsilon)}, ~n' \neq n} x_{n'} \right| & \leq  \sum_{n' \equiv n \textrm{~mod~} t_i \cdot s_{j'},~n' \notin S^{\rm opt}_{(k / \epsilon)},~n' \neq n} y_{n'} \nonumber \\ & \leq \sum_{n' \equiv n \textrm{~mod~} s_{j'},~n' \notin S^{\rm opt}_{(k / \epsilon)}, ~n' \neq n} y_{n'} = \left| \left( \mathcal{M}_{s_1,K,n} \cdot \vec{y} \right)_{j'} - y_{n} \right| \leq \bar{\delta}.
\nonumber
\end{align}
Finally, to verify property $(3)$ we can bound $\left| \left( r_{j',n \textrm{~mod~} s_{j'}} \circledast \tilde{r}_{i, h} \right) \cdot \vec{x} ~\right|$ from above for all $i \in [1,\lambda] \cap \mathbbm{N}$ and $h \in [0,t_i) \cap \left( \mathbbm{N} - \{ n \textrm{~mod~} t_i \} \right)$ by 
$$\left| \sum_{n' \equiv n \textrm{~mod~} s_{j'},~n' \equiv h \textrm{~mod~} t_{i},~n' \notin S^{\rm opt}_{(k / \epsilon)}} x_{n'} \right| \leq \sum_{n' \equiv n \textrm{~mod~} s_{j'},~n' \notin S^{\rm opt}_{(k / \epsilon)}, ~n' \neq n} y_{n'} = \left| \left( \mathcal{M}_{s_1,K,n} \cdot \vec{y} \right)_{j'} - y_{n} \right| \leq \bar{\delta}.$$
Hence, we can see that all three properties will indeed hold for at least $\frac{K}{2}$ rows of $\mathcal{M}_{s_1,K,n}$.

\section{Proof of Theorem~\ref{thm:Recon2}}
\label{app:proof_thm:Recon2}

Let $\delta$ be defined as 
$$\delta = \frac{ \epsilon \cdot \left\| \vec{\hat{f}} - \vec{\hat{f}}^{\rm ~opt}_{(k/\epsilon)} \right\|_1}{k} + \left\| ~\hat{f} - \bar{\hat{f}} ~ \right\|_1.$$
Furthermore, suppose $j \in [1,K] \cap \mathbbm{N}$ and $h \in [0,s_j)$ correspond to an $\omega_{j,h} \in \left( - \left\lceil \frac{N}{2} \right\rceil, \left\lfloor \frac{N}{2} \right\rfloor \right] \cap \mathbbm{Z}$ which is reconstructed more than $\frac{K}{2}$ times by line 12 of Algorithm~\ref{alg:reconstruct2}.  As a consequence of Lemmas~\ref{lem:MultError} and~\ref{lem:Rprops} we can see than more than half of the entries of $\mathcal{G}_{\lambda,K,\omega_{j,h}} \tilde{\Psi} \vec{A}$ produced in line 4 will satisfy 
$\left| \left(\mathcal{G}_{\lambda,K,\omega_{j,h}} \tilde{\Psi} \vec{A} \right)_j - \hat{f}_{\omega_{j,h}} \right| \leq \delta$.  Therefore, the $x_{\omega_{j,h}}$ value produced by lines 16 and 17 will have 
\begin{equation}
\left| x_{\omega_{j,h}} - \hat{f}_{\omega_{j,h}} \right| \leq \sqrt{2} \cdot \delta.
\label{eqn:ApproxE2}
\end{equation}

Since Equation~\ref{eqn:ApproxE2} will hold for all $\omega \in \left( - \left\lceil \frac{N}{2} \right\rceil, \left\lfloor \frac{N}{2} \right\rfloor \right] \cap \mathbbm{Z}$ reconstructed more than $\frac{K}{2}$ times, we can begin to bound the approximation error by 
\begin{align}
\left\| \vec{\hat{f}} - \vec{x}_{S} \right\|_2 &~\leq~ \left\| \vec{\hat{f}} - \vec{\hat{f}}_{S} \right\|_2 + \left\| \vec{\hat{f}}_{S} - \vec{x}_{S} \right\|_2 ~\leq~ \left\| \vec{\hat{f}} - \vec{\hat{f}}_{S} \right\|_2 + 2\sqrt{k} \cdot \delta \nonumber \\
&~=~\sqrt{\left\| \vec{\hat{f}} - \vec{\hat{f}}^{\rm~ opt}_k \right\|^2_2 + \sum_{\omega \in S^{\rm opt}_k - S} \left| \hat{f}_{\omega} \right|^2 - \sum_{\tilde{\omega} \in S - S^{\rm opt}_k } \left| \hat{f}_{\tilde{\omega}} \right|^2} + 2\sqrt{k} \cdot \delta.
\label{eqn:ApproxEE2}
\end{align}
In order to make additional progress on Equation~\ref{eqn:ApproxEE2} we must now consider the possible magnitudes of $\vec{\hat{f}}$ entries at indices in $S - S^{\rm opt}_k$ and $S^{\rm opt}_k - S$.

Suppose $\omega \in S^{\rm opt}_k - S \neq \emptyset$.  In this case either $(i)$ $\left| \hat{f}_{\omega} \right| \leq 4 \delta$, or $(ii)$ $\left| \hat{f}_{\omega} \right| > 4 \delta$ in which case Lemma~\ref{lem:FreqIdent} guarantees that $\omega$ will be identified by lines 6 through 14 of Algorithm~\ref{alg:reconstruct2}.  Once identified, an $\bar{\omega} \in S^{\rm opt}_k$ will always be placed in $S$ unless at least $k+1$ other distinct identified elements, $\tilde{\omega} \notin S^{\rm opt}_k$, have the property that $|x_{\tilde{\omega}}| \geq |x_{\bar{\omega}}|$.  Thus, if $(ii)$ occurs then 
\begin{equation}
\left|\hat{f}_{\omega_{k}} \right| + \sqrt{2} \cdot \delta ~\geq~ \left|\hat{f}_{\tilde{\omega}} \right| + \sqrt{2} \cdot \delta ~\geq~ \left|\hat{f}_{\omega} \right| - \sqrt{2} \cdot \delta ~\geq~ \left|\hat{f}_{\omega_{k}} \right| - \sqrt{2} \cdot \delta
\label{eqn:ApproxE3}
\end{equation}
will hold for all $\tilde{\omega} \in S - S^{\rm opt}_k$.
The end result is that if $\omega \in S^{\rm opt}_k - S$ then either $\left|\hat{f}_{\omega} \right| \leq 4 \delta$, or else $\hat{f}_{\omega}$ is roughly the same magnitude as $\hat{f}_{\omega_{k}}$ (up to a $O(\delta)$ tolerance).  Furthermore, because line 20 chooses $2k$ elements for $S$ whenever possible, we can see that
$S - S^{\rm opt}_k$ must contain at least 
$$2 \cdot \left| \left(S^{\rm opt}_k \cap \left\{ \omega ~\Big|~ \left|\hat{f}_{\omega} \right| > 4 \delta \right\} \right) - S \right|$$ 
elements, $\tilde{\omega}$, all of which satisfy Equation~\ref{eqn:ApproxE3} for every $\omega \in \left( S^{\rm opt}_k \cap \left\{ \omega ~\Big|~ \left|\hat{f}_{\omega} \right| > 4 \delta \right\} \right) - S$.  We are now ready to give Equation~\ref{eqn:ApproxEE2} further consideration.

If $S^{\rm opt}_k - S = \emptyset$ we are finished.  Otherwise, if $S^{\rm opt}_k - S \neq \emptyset$, we can bound the squared $l^2$-norm of $\vec{\hat{f}}_{S - S^{\rm opt}_k}$ from below by
$$\sum_{\tilde{\omega} \in S - S^{\rm opt}_k } \left| \hat{f}_{\tilde{\omega}} \right|^2 \geq 2 \cdot \left| \left(S^{\rm opt}_k \cap \left\{ \omega ~\Big|~ \left|\hat{f}_{\omega} \right| > 4 \delta \right\} \right) - S \right| \cdot \left( \left| \hat{f}_{\omega_{k}} \right| - 2 \sqrt{2} \cdot \delta \right)^2~=~{\bf A}.$$
Furthermore, we can upper bound the squared $l^2$-norm of $\vec{\hat{f}}_{S^{\rm opt}_k - S}$ by
$$\left|\left(S^{\rm opt}_k \cap \left\{ \omega ~\Big|~ \left|\hat{f}_{\omega} \right| > 4 \delta \right\} \right) - S \right| \cdot \left( \left| \hat{f}_{\omega_{k}} \right| + 2 \sqrt{2} \cdot \delta \right)^2 + \left|\left(S^{\rm opt}_k \cap \left\{ \omega ~\Big|~ \left|\hat{f}_{\omega} \right| \leq 4 \delta \right\} \right) - S \right| \cdot 16 \delta^2 \geq \sum_{\omega \in S^{\rm opt}_k - S} \left| \hat{f}_{\omega} \right|^2.$$
Let ${\bf B} = \left|\left(S^{\rm opt}_k \cap \left\{ \omega ~\Big|~ \left|\hat{f}_{\omega} \right| > 4 \delta \right\} \right) - S \right| \cdot \left( \left| \hat{f}_{\omega_{k}} \right| + 2 \sqrt{2} \cdot \delta \right)^2$.  We will now concentrate on bounding $${\bf C} = \sum_{\omega \in S^{\rm opt}_k - S} \left| \hat{f}_{\omega} \right|^2 - \sum_{\tilde{\omega} \in S - S^{\rm opt}_k } \left| \hat{f}_{\tilde{\omega}} \right|^2.$$

If ${\bf A} \geq {\bf B}$ then ${\bf C} \leq 16 k \cdot \delta^2$.  Otherwise, if ${\bf A} < {\bf B}$ then 
$$\left| \hat{f}_{\omega_{k}} \right|^2 - 12 \sqrt{2} \delta \cdot \left| \hat{f}_{\omega_{k}} \right| + 8 \delta^2 < 0$$ 
which can only happen if $\left| \hat{f}_{\omega_{k}} \right| \in \left( (6 \sqrt{2} - 8) \cdot \delta, (6 \sqrt{2} + 8) \cdot \delta \right)$.  Hence, ${\bf A} < {\bf B}$ implies that ${\bf C} \leq k \cdot \left(8 \sqrt{2} + 8 \right)^2 \cdot \delta^2$.

Finishing our error analysis, we can see that in the worst possible case Equation~\ref{eqn:ApproxEE2} will remain bounded by
$$\left\| \vec{\hat{f}} - \vec{x}_{S} \right\|_2 \leq \sqrt{\left\| \vec{\hat{f}} - \vec{\hat{f}}^{\rm~ opt}_k \right\|^2_2 + k \cdot \left(8 \sqrt{2} + 8 \right)^2 \cdot \delta^2} ~+~ 2\sqrt{k} \cdot \delta ~\leq~ \left\| \vec{\hat{f}} - \vec{\hat{f}}^{\rm~ opt}_k \right\|_2 ~+~ 22\sqrt{k} \cdot \delta.$$
The error bound stated in Equation~\ref{eqn:Recon2Error} follows.  The upper bound on the number of point evaluations of $f$ follows from an application of Lemma~\ref{lem:RowBound} and Theorem~\ref{thm:RowBound} with $c = 4$.    

We will begin bounding the runtime of Algorithm~\ref{alg:reconstruct2} by bounding the runtime of lines 15 through 21.  Line 12 of Algorithm~\ref{alg:reconstruct2} will be executed a total of $O \left( \frac{k^2 \cdot \log^2_{(k/\epsilon)} N \cdot \log \left(\frac{k \cdot \ln N}{\epsilon} \right)}{\epsilon^2} \right)$ times (see Equation~\ref{eqn:NumRows}).  Therefore, lines 16 and 17 will be executed $O \left( \frac{k \cdot \log_{(k/\epsilon)} N \cdot \log \left(\frac{k \cdot \ln N}{\epsilon} \right)}{\epsilon} \right)$ times apiece.  Each such median operation can be accomplished in $O(K \cdot \lambda)$ time using a median-of-medians algorithm (e.g., see \cite{MedianMedian}).  Therefore, the total runtime of lines 15 through 21 will be $O \left( \frac{k^2 \cdot \log^2_{(k/\epsilon)} N \cdot \log \left(\frac{k \cdot \ln N}{\epsilon} \right) \cdot \log(\epsilon N/k)}{\epsilon^2 \cdot \log \log (\epsilon N/k)} \right)$.  Turning our attention to lines 6 through 14, we note that their runtime will be dominated by the $O \left( \frac{k^2 \cdot \log^2_{(k/\epsilon)} N \cdot \log \left(\frac{k \cdot \ln N}{\epsilon} \right) \cdot \log(\epsilon N/k)}{\epsilon^2 \cdot \log \log (\epsilon N/k)} \right)$ executions of line 9.  Therefore, the total runtime of lines 6 through 14 will be $O \left( \frac{k^2 \cdot \log^2_{(k/\epsilon)} N \cdot \log \left(\frac{k \cdot \ln N}{\epsilon} \right) \cdot \log^2(\epsilon N/k)}{\epsilon^2 \cdot \log \log (\epsilon N/k)} \right)$ (see \cite{IDSFA} and Lemma~\ref{lem:RowBound}).  The stated overall runtime of lines 6 through 21 follows.

\section{Proof of Corollary~\ref{cor:Recon2}}
\label{app:proof_cor:Recon2}

Define $\vec{\widehat{|f|}} \in \mathbbm{R}^N$ by
$$\left(\vec{\widehat{|f|}}\right)_{\omega} = \left| \hat{f}_{\omega} \right|~\textrm{for all}~\omega \in \left( - \left\lceil \frac{N}{2} \right\rceil, \left\lfloor \frac{N}{2} \right\rfloor \right] \cap \mathbbm{Z}.$$
Clearly $\vec{\hat{f}}$ and $\vec{\widehat{|f|}}$ will both have the same optimal $(k/ \epsilon)$-term support subset, $S^{\rm opt}_{(k/ \epsilon)} \subset [0,N) \cap \mathbbm{N}$.  Similarly, it is easy to see that $\left\| \vec{\hat{f}} - \vec{\hat{f}}^{\rm opt}_{(k/\epsilon)} \right\|_1 = \left\| \vec{\widehat{|f|}} - \left(\vec{\widehat{|f|}} \right)_{(k/\epsilon)}^{\rm opt} \right\|_1$.  Apply Corollary~\ref{cor:RandSj} with $c = 14$, $\vec{x} = \vec{\widehat{|f|}}$, and $S = \left( - \left\lceil \frac{N}{2} \right\rceil, \left\lfloor \frac{N}{2} \right\rfloor \right] \cap \mathbbm{Z}$ to obtain $\tilde{S}$, a multiset of $\left\lceil 21 \cdot \ln \left( \frac{N}{1 - \sigma} \right) \right\rceil$ $s_j$ values.  With probability at least $\sigma$ more than half (with multiplicity) of the entries of $\mathcal{M}_{\tilde{S}},\omega \cdot \vec{\widehat{|f|}}$ will estimate the $\omega^{\rm th}$ entry of $\vec{\widehat{|f|}}$ to within $(\epsilon/k) \cdot \left\| \vec{\hat{f}} - \vec{\hat{f}}^{\rm ~opt}_{(k/\epsilon)} \right\|_1$ precision for all $\omega \in S$. 

Given the last paragraph, it is not difficult to see that with probability at least $\sigma$ a result analogous to that of Lemma~\ref{lem:Rprops} will hold for $\mathcal{R}_{\lambda,\tilde{S}} \cdot \vec{\hat{f}}$.  That is, with probability at least $\sigma$ the following will hold for all $\omega \in S$:  The majority (when counted with multiplicity) of $\mathcal{M}_{\tilde{S}},\omega$ rows, $\vec{r} \in \{ 0, 1\}^N$, will have $\left( \vec{r} \circledast \vec{s} \right) \cdot \vec{\hat{f}} \approx \hat{f}_{\omega}$ for a given row, $\vec{s}$, of $\mathcal{N}_{\lambda,s_1}$ if and only if $\vec{s}$ is also a row of $\mathcal{N}_{\lambda,s_1,\omega}$.  Furthermore, $\mathcal{R}_{\lambda,\tilde{S}} \cdot \vec{\hat{f}}$ can still be approximately computed using Algorithm~\ref{alg:fastMult} if only the unique $s_j$ values in $\tilde{S}$ are given as relatively prime $s_j$-inputs.  In this case a result analogous to Lemma~\ref{lem:MultError} will also still hold since we will merely be computing a subset of the previously calculated vector entries.  Finally, by inspecting the proof of Lemma~\ref{lem:FreqIdent} we can see that an almost identical result (with $K$ replaced by the $l$ value from Corollary~\ref{cor:RandSj}) will hold any time $\mathcal{R}_{\lambda,\tilde{S}} \cdot \vec{\hat{f}}$ satisfies the aforementioned variants of both Lemmas~\ref{lem:Rprops} and~\ref{lem:MultError}.

Taken all together, we can see that with probability at least $\sigma$ both of the following statements will be true:  First, all at most $N$ $x_{\omega_{j,h}}$ values ever produced by lines 16 and 17 of Algorithm~\ref{alg:reconstruct2} will have 
$$\left| x_{\omega_{j,h}} - \hat{f}_{\omega_{j,h}} \right| \leq \sqrt{2} \cdot \left(\frac{ \epsilon \cdot \left\| \vec{\hat{f}} - \vec{\hat{f}}^{\rm ~opt}_{(k/\epsilon)} \right\|_1}{k} + \left\| ~\hat{f} - \bar{\hat{f}} ~ \right\|_1 \right).$$
Second, a variant of Lemma~\ref{lem:FreqIdent} will ensure that all $\omega \in S$ with
$$\left| \hat{f}_{\omega} \right| ~>~ 4 \cdot \left( \frac{ \epsilon \cdot \left\| \vec{\hat{f}} - \vec{\hat{f}}^{\rm ~opt}_{(k/\epsilon)} \right\|_1}{k} + \left\| ~\hat{f} - \bar{\hat{f}} ~ \right\|_1 \right)$$
are reconstructed by lines 6 through 14 of Algorithm~\ref{alg:reconstruct2} more than $\frac{\left\lceil 21 \cdot \ln \left( \frac{N}{1 - \sigma} \right) \right\rceil}{2}$ times.  The Equation~\ref{eqn:Recon2Error} error bound now follows from the proof of Theorem~\ref{thm:Recon2}.

We upper bound the number of required function evaluations by bounding the number of rows for a particular $\mathcal{R}_{\lambda,\tilde{S}}$ matrix constructed with $\left\lceil 21 \cdot \ln \left( \frac{N}{1 - \sigma} \right) \right\rceil$ randomly chosen $s_j$ values defined as in Equations~\ref{eqn:Primes} --~\ref{eqn:Def_k} \big(with $K = 14 \cdot (k / \epsilon) \left \lfloor \log_{s_1} N \right \rfloor + 1$ \big).  In this case Equation~\ref{eqn:q+Kbound} together with results from \cite{PrimeBound} tell us that $s_K$ is itself bounded above by Equation~\ref{eqn:SKboundApp}.  The final upper bound on the number of point evaluations of $f$ then follows from an application of Lemma~\ref{lem:RowBound}.  Note that the product of the Lemma~\ref{lem:RowBound} row bound with $\left\lceil 21 \cdot \ln \left( \frac{N}{1 - \sigma} \right) \right\rceil$ and  Equation~\ref{eqn:SKboundApp} provides a concrete upper bound for the number of point evaluations of $f$.

We will begin bounding the runtime of Algorithm~\ref{alg:reconstruct2} by bounding the runtime of lines 15 through 21.  Line 12 of Algorithm~\ref{alg:reconstruct2} will be executed a total of $O \left( \frac{k \cdot \log \left( \frac{N}{1 - \sigma} \right) \cdot \log_{(k/\epsilon)} N \cdot \log \left(\frac{k \cdot \log N}{\epsilon} \right)}{\epsilon} \right)$ times (see Equation~\ref{eqn:SKboundApp}).  Therefore, lines 16 and 17 will be executed $O \left( \frac{k \cdot \log_{(k/\epsilon)} N \cdot \log \left(\frac{k \cdot \ln N}{\epsilon} \right)}{\epsilon} \right)$ times apiece.  Each such median operation can be accomplished in $O\left( \log \left( \frac{N}{1 - \sigma} \right) \cdot \lambda \right)$ time using a median-of-medians algorithm (e.g., see \cite{MedianMedian}).  Therefore, the total runtime of lines 15 through 21 will be $O \left( \frac{k \cdot \log \left( \frac{N}{1 - \sigma} \right) \cdot \log_{(k/\epsilon)} N \cdot \log \left(\frac{k \cdot \log N}{\epsilon} \right) \cdot \log(\epsilon N/k)}{\epsilon \cdot \log \log (\epsilon N/k)} \right)$.  Turning our attention to lines 6 through 14, we note that their runtime will be dominated by the $O \left( \frac{k \cdot \log \left( \frac{N}{1 - \sigma} \right) \cdot \log_{(k/\epsilon)} N \cdot \log \left(\frac{k \cdot \log N}{\epsilon} \right) \cdot \log(\epsilon N/k)}{\epsilon \cdot \log \log (\epsilon N/k)} \right)$ executions of line 9.  Therefore, the total runtime of lines 6 through 14 will be $O \left( \frac{ k \cdot \log \left( \frac{N}{1 - \sigma} \right) \cdot \log_{(k/\epsilon)} N \cdot \log \left(\frac{k \cdot \log N}{\epsilon} \right) \cdot \log^2(\epsilon N/k)}{\epsilon \cdot \log \log (\epsilon N/k)} \right)$ (see Lemma~\ref{lem:RowBound}).  The stated overall runtime of lines 6 through 21 follows.

\section{Proof of Theorem~\ref{thm:MultDimRecov} }
\label{app:MultDimRecov}

Suppose we want to resolve at least $M$ frequencies in each of $D$ dimensions (i.e., we want to approximate the $D$-dimensional array $\vec{\hat{f}} \in \mathbbm{C}^{M^D}$).  We begin by choosing the smallest $\tilde{D} \in \mathbbm{N}$ such that  
$$\prod^{\tilde{D}-D+1}_{j=1} p_j ~>~ \left(MD \right)^D.$$
The first paragraph of Appendix~\ref{app:proof_lem:RowBound} reveals that 
\begin{equation}
\tilde{D} ~<~  3 \cdot \frac{D \cdot \ln( MD )}{\ln (D \cdot \ln( MD ))} + D ~=~ \frac{D \cdot \ln( MD )}{\ln (D \cdot \ln( MD ))} \cdot 3 \left( 1 + \frac{ \ln D + \ln \ln( MD )}{\ln( MD )} \right) ~=~ O\left( \frac{D \cdot \log( MD )}{\log (D \cdot \log( MD ))}\right).
\label{eqn:MULT_D_bound}
\end{equation}
Furthermore, we can see from \cite{PrimeBound} that

\begin{equation}
\ln \left( \prod^{\tilde{D}}_{j=1} p_j \right) ~\leq~ D \cdot \ln( MD ) + \sum^{\tilde{D}}_{j=\tilde{D}-D+1} \ln p_j ~\leq~ D \cdot \left( \ln( MD ) + \ln p_{\tilde{D}} \right) ~\leq~D \cdot \left( \ln( MD ) + \ln \left( \tilde{D} \cdot ( \ln \tilde{D} + \ln \ln \tilde{D} \right) \right)
\label{eqn:MULT_N_bound}
\end{equation}
for $D \geq 2$ and $(MD)^D \geq 2,310$.  Thus, $\log \left( \prod^{\tilde{D}}_{j=1} p_j \right)$ is generally $O(D \cdot \log ( M \cdot D))$.  We will use these first $\tilde{D}$ primes to help define our new one-dimensional function $f_{\rm new}$ (see Equation~\ref{def:fnew} above).

Set $\tilde{D}_0 = 0$ and recursively define $\tilde{D}_d$ to be such that
$$\prod^{\tilde{D}_d - 1}_{j = \tilde{D}_{d-1}+1} p_j ~\leq~ MD ~<~ \prod^{\tilde{D}_d}_{j = \tilde{D}_{d-1}+1} p_j$$
for all $1 \leq d \leq D$.  We then define the $D$ pairwise relatively prime values required for approximating each Fourier coefficient, $\vec{\hat{f}}\left(g^{-1}(\omega) \right) \in \mathbbm{C}$, via Equation~\ref{eqn:fnewFourier} to be 
$$P_d = \prod^{\tilde{D}_d}_{j = \tilde{D}_{d-1}+1} p_j$$
for $1 \leq d \leq D$.  Set $\tilde{N} ~=~ \prod^{D}_{d=1} P_d ~\leq~ \prod^{\tilde{D}}_{j=1} p_j$.  The stated runtime and sampling bounds follow.  

We are now in the position to apply any of Theorem~\ref{thm:Recon1}, Corollary~\ref{cor:RandRecon1}, Theorem~\ref{thm:Recon2}, or Corollary~\ref{cor:Recon2} to approximate $\vec{\hat{f}}_{\rm new} \in \mathbbm{C}^{\tilde{N}}$.  Upon applying any of these four results to $f_{\rm new}$ we will obtain (either deterministically, or randomly with high probability) a $2k$-sparse $\vec{x}_{S} \in \mathbbm{C}^{\tilde{N}}$ satisfying
$$\left\| \vec{\hat{f}}_{\rm new} - \vec{x}_{S} \right\|_2 \leq \left\| \vec{\hat{f}}_{\rm new} - \vec{\hat{f}}^{\rm~ opt}_{{\rm new}~k} \right\|_2 + \frac{22\epsilon \cdot \left\| \vec{\hat{f}}_{\rm new} - \vec{\hat{f}}^{\rm~ opt}_{{\rm new}~(k/\epsilon)} \right\|_1}{\sqrt{k}} + 22\sqrt{k} \cdot \left\| ~\hat{f}_{\rm new} - \bar{\hat{f}}_{\rm new} ~ \right\|_1.$$
Recall that we have only guaranteed that $\hat{f}$'s Fourier coefficients for $\left( [-M/2, M/2] \cap \mathbbm{Z} \right)^D$ map into $\hat{f}_{\rm new}$'s Fourier coefficients for $[-\tilde{N}/2, \tilde{N}/2] \cap \mathbbm{Z}$ (although many others will as well).  Thus, for simplicity, we assumed that $f$ is bandlimited when translating these error bounds back into terms of $\hat{f}$.  Given this assumption we have $\left\| ~\hat{f}_{\rm new} - \bar{\hat{f}}_{\rm new} ~ \right\|_1 = 0$.  The stated error bound now follows from the fact that $g$ is a bijection.  

\end{document}